\documentclass{amsart}
  \usepackage{url}
  \usepackage{amssymb}
  \input xypic
  \xyoption{all}

  \usepackage{amsmath, amsthm,  amsfonts}

  \theoremstyle{plain}
  \newtheorem{theorem}{Theorem}[section]
  \newtheorem{lemma}[theorem]{Lemma}
  \newtheorem{proposition}[theorem]{Proposition}
  
  \newtheorem{corollary}{Corollary}[section]
  \newtheorem{definition}[theorem]{Definition}
  \newtheorem{remark}[theorem]{Remark}
  
  \newtheorem{note}{Note}[section]

  \numberwithin{equation}{section}
  \numberwithin{figure}{section}

  \renewcommand{\cH}{{\mathcal H}}
  \newcommand{\cB}{{\mathcal B}}
  \newcommand{\cK}{{\mathcal K}}
  
  \newcommand{\cA}{{\mathcal A}}
  \newcommand{\cE}{{\mathcal E}}

  \newcommand{\cP}{{\mathcal P}}

\newcommand{\cFr}{{\mathcal Fr }}

   \newcommand{\ba}{\begin{eqnarray}}
   \newcommand{\na}{\end{eqnarray}}
   \newcommand{\ban}{\begin{eqnarray*}}
   \newcommand{\nan}{\end{eqnarray*}}

  \newcommand{\CC}{{\mathbb C}}
  \newcommand{\RR}{{\mathbb R}}
  \newcommand{\ZZ}{{\mathbb Z}}
  
  \renewcommand{\AA}{{\mathbb A}}

  \newcommand{\PP}{{\mathbb P}}

    \newcommand{\disp}{\displaystyle}

  \def\cancel#1#2{\ooalign{$\hfil#1\mkern1mu/\hfil$\crcr$#1#2$}}
\def\Dirac{\mathpalette\cancel D}

\begin{document}

\title[Thom isomorphism and Push-forward map  in  twisted K-theory]
{Thom isomorphism and Push-forward map  in twisted K-theory}

  \author[A.L. Carey]{Alan L. Carey}
\address[Alan L. Carey]
  {Mathematical Sciences Institute\\
  Australian National University\\
  Canberra ACT 0200 \\
  Australia}
  \email{acarey@maths.anu.edu.au}

  \author[Bai-Ling Wang]{Bai-Ling Wang}
  \address[Bai-Ling Wang]
  {Mathematical Sciences Institute\\
  Australian National University\\
  Canberra ACT 0200 \\
  Australia}
  \email{wangb@maths.anu.edu.au}

  \thanks{This work is  supported by  the Australian
  Research Council Discovery Project DP0449470.
   The authors like to thank Ulrich Bunke,  David Evans, Michael Murray and Danny Stevenson for
   many useful conversations.
   }

  \begin{abstract} We  establish the
  Thom isomorphism in twisted K-theory for any  real vector bundle
  and develop the push-forward map  in twisted K-theory for
  any differentiable map $f: X\to Y$ (not necessarily K-oriented).
  The push-forward map  generalizes  the push-forward map  in  ordinary
  K-theory for   any $K$-oriented  map  and the Atiyah-Singer
  index theorem of Dirac operators on Clifford modules. For $D$-branes satisfying Freed-Witten's
  anomaly cancellation condition in a manifold with a non-trivial
  $B$-field, we associate a canonical element in the twisted
  K-group to get the so-called D-brane charges. 
  \end{abstract}

\maketitle

\tableofcontents

  \section{Introduction}

  In complex K-theory, the push-forward map $f_!^{c_1}: K^*(X) \to K^*(Y)$ was established
  by Atiyah-Hirzebruch \cite{AH} for any differentiable   $c_1$-map $f: X\to Y$
  where  $c_1\in H^2(X, \ZZ)$ satisfies
  \[
 c_1  \equiv w_2(X) -f^*w_2(Y) \ (mod  \ 2)
  \]
  with $w_2(X)$ and  $w_2(Y)$ being the second Stiefel-Whitney classes of $X$ and $Y$ respectively. 
  Such a map $f: X\to Y$ is called K-oriented.

  The  push-forward map $f_{!}^{c_1}$, also called the Gysin homomorphism in K-theory, gives rise to
  the following general Riemann-Roch theorem, for $a\in K^*(X)$,
  \[
  Ch\bigl(f_{!}^{c_1}(a)\bigr) \hat{A}(X) =  f_*\bigl(Ch(a)e^{\frac{c_1(f)}{2}}\hat{A}(Y)\bigr),
  \]
  where $f_*$ is the Gysin homomorphism in ordinary cohomology theory,
  $Ch$ is the Chern character for K-theory, $\hat{A}(X)$ and $\hat{A}(Y)$ are the
  A-hat genus (which can be expressed in terms of the Pontrjagin classes) of $X$ and $Y$ respectively.
  The set of K-orientations for $f: X\to Y$ is an affine space modeled on $2H^2(X, \ZZ)$, hence
  $H^2(X, \ZZ)$ acts (and  transitively if $H^2(X, \ZZ)$ has no 2-torsion elements)
  via
  \[
  c_1 \mapsto c_1 + 2c,
  \]
  for $c\in H^2(X, \ZZ)$.
  The dependence of  the push-forward map $f_{!}^{c_1}$ on the K-orientation $c_1$ can be described
  by
  \[
  f_{!}^{c_1+ 2c } (a)
 =    f_{!}^{c_1} (a\cdot [L_c])
 \]
 where $[L_c] \in K^0(X)$ is the K-element defined by the equivalence class of the line
 bundle over $X$ with first Chern class $c\in  H^2(X, \ZZ)$.

  In this paper, we will develop the push-forward map in twisted K-theory for any
  differentiable  map $f: X\to Y$, which is not necessarily K-oriented. 
  We aim to present the results in a fashion which demonstrates that for
  general manifolds and maps twisted K-theory provides the natural framework 
  in which to study the push-forward map,   Thom isomorphism and topological index.
  
  Recall that  twisted
  K-theory for a smooth manifold $X$ with a class $\sigma \in H^3(X, \ZZ)$
  is defined using  a locally trivial projective Hilbert
  bundle over $X$ with infinite dimensional separable fibers, $P_\sigma$, whose
   Dixmier-Douady class 
  is $\sigma \in H^3(X, \ZZ)$:
  \[
  K^*_{P_\sigma }(X) =K^0_{P_\sigma }(X)   \oplus K^1_{P_\sigma}(X).
  \]
  As a $\ZZ_2$-graded cohomology theory it was fully developed by Atiyah and Segal
  in \cite{AS}. Earlier approaches may be found in  \cite{DK}\cite{Ros}\cite{BCMMS}.
  Automorphisms of the $PU(\cH)$-principal bundle $ \cP_\sigma$ 
  result in a natural action of the Piccard group  $Pic(X)$ 
  on the twisted K-group $K^*_{\sigma }(X)$. 
   Up to this natural action of $Pic(X)$, $K^*_{P_\sigma }(X)$ depends
  only on the Dixmier-Douady class of $P_\sigma$. With this understood, we often
  denote the twisted  K-group by $K^*_\sigma(X)$. Though in practice, we always 
  fix a $PU(\cH)$-principal bundle $ \cP_\sigma$ with  Dixmier-Douady class $\sigma$.

  Twisted K-theory for manifolds carrying a projective Hilbert
  bundle is a homotopy invariant and satisfies the Mayer-Vietoris property.  
  Twisted K-theory  has attracted  considerable interest.  
  It is intimately related to D-brane physics in superstring theory.
  Moreover the result  of Freed-Hopkins-Teleman 
  \cite{FHT} which identifies the equivariant twisted K-theory of a compact Lie
  group $G$ with the Verlinde ring of projective representations of the loop group $LG$
  has focussed the attention of many researchers.

 We now recall the viewpoint of \cite{BCMMS} which we will exploit later in the paper.
  Denote by
  $W_3(f)$ the image of $w_2(X)- f^*(w_2(Y)) \in H^2(X, \ZZ_2)$ under the Bockstein homomorphism
  \[
  \beta: H^2(X, \ZZ_2) \longrightarrow H^3(X, \ZZ).
  \]
  For $f:X \to pt$, $W_3(f)$ is the third integral Stiefel-Whitney class $W_3(X)$ of $X$.
 Now  \cite{BCMMS} regards
  $K^0_{ W_3(X)}(X)$ as the Grothendieck group formed from  bundle gerbe modules
  of a bundle gerbe over $X$ with Dixmier-Douady class $W_3(X)$.  
   
   The main theorem of this paper is   the following.

 \vspace{4mm}

  \begin{theorem}  Let  $f: X\to Y$  be a  differentiable   map and $\sigma
  \in H^3(Y,\ZZ)$.
  There exists a push-forward map,  compatible with the natural  action of the
  Piccard group $Pic (Y)$, 
  \[
  f_{!}: K^*_{f^*\sigma + W_3(f)}(X) \longrightarrow K^*_{\sigma }(Y)
  \]
  with grading shifted by $dim X - dim Y (mod \ 2)$. In particular,
  \begin{enumerate}
  \item if $f: X\to Y$ is K-oriented, i.e., $W_3(f) =0$,  with $c_1\in H^2(X, \ZZ)$ satisfying
  \[
  c_1  \equiv w_2(X) -f^*w_2(Y) \ (mod  \ 2),
  \]
  there exists a natural push-forward map
  \[
  f_{!}^{c_1}: K^*_{f^*\sigma }(X) \longrightarrow K^*_{\sigma }(Y)
  \]
  which agrees with the push-forward map for ordinary K-theory when $\sigma$ is trivial.
  \item if $f: X \to pt$, then
  \[
  f_{!}: K^*_{  W_3(X)}(X) \longrightarrow K^*(pt) \cong \ZZ
  \]
  is the Atiyah-Singer  index theorem for Dirac operators on Clifford modules
  (Theorem 4.3 in \cite{BGV}),
 see also
 \cite{MurSin} and the introduction of
  \cite{AS}.
  \item if $\iota: Q \to X$ is a closed submanifold $Q$ of  a spin manifold $X$ such that
  \[\iota^*\sigma +W_3(Q) =0,
  \]
  then the push-forward map defines the charge of a $D$-brane   $(Q, \xi)$ supported on $Q$
  for $\xi\in K(Q)$
  \[
  \iota_!: K(Q)\longrightarrow K^*_{  \sigma} (X).
  \]
  \item if  $f: X\to Y$  is a fibration with a closed oriented fiber, then the topological 
  index in twisted K-theory is given by 
  \ba\label{topo:tw-index}
  \xymatrix{
  K^*_{\pi^*f^*\sigma} (T(X/Y)) \ar[r]^{\cong} 
  &  K^*_{f^*\sigma + W_3(f)}(X) \ar[r]^{\ \ \ f_!} & K^*_{\sigma }(Y),
  }
  \na
   where $\pi: T(X/Y) \to X$ is the vertical tangent bundle of $f: X\to Y$ and $W_3(\pi) = W_3(f)$. 
    \end{enumerate}
    \end{theorem}

For a differentiable   map $f: X \to Y$ with the $K$-orientation condition (i.e., $W_3(f) =0$),  
the above push-forward map in twisted K-theory 
is also established in \cite{MMSin}\cite{BEM} by using a $C^*$-algebraic approach. The topological 
index (\ref{topo:tw-index})  in twisted K-theory when $W_3(f) =0$ is developed in 
 \cite{MMSin} using a different method.

One of the main technical issues is to establish the Thom isomorphism in twisted
K-theory  for  general twistings $\sigma \in H^3(M, \ZZ)$,
 which is done in section \ref{Thom}. This also generalizes
the earlier result of Donovan-Karoubi in \cite{DK} for a torsion twisting. 

In order to give a uniform treatment in both the torsion and non-torsion cases
we reprove the Thom isomorphism in twisted
K-theory  for torsion twisting using the language of bundle gerbes and bundle gerbe modules.
For non-torsion twisting $\sigma$, let $\cP_\sigma$ be the corresponding principal $PU(\cH)$-bundle.
We associate a bundle $$\cA_\sigma= \cP_\sigma\times_{PU(\cH)} End_{HS}(\cH)$$ 
of Hilbert algebras, fiberwisely
consisting of the Hilbert algebra  $End_{HS}(\cH)$ of Hilbert-Schmidt operators with the conjugation 
action of $PU(\cH)$. This  bundle of Hilbert algebras plays the role that the
 Azumaya bundle plays for the torsion case. 

Let $U_2$ be the subgroup of $U(\cH)$ of unitaries of the form  $1+$ Hilbert-Schmidt.
 The complex Hilbert  bundles with structure group $U_2$ and admitting a fiberwise $\cA_\sigma$-action  
  form  an additive category $\cE_{U_2}^{\cA_\sigma}(M)$, whose
 Grothendieck group is isomorphic to the twisted K-group $K_\sigma(M)$.
  From this we show the existence  of the Thom isomorphism in twisted K-theory.

\vspace{4mm}

\begin{theorem}
Suppose that $\pi: V\to M$ is an oriented real
vector bundle over $M$ of even rank with positive definite
quadratic form.  There is a natural isomorphism
$$
K^0_{\sigma+W_3(V)}(M) \cong K^0_{\pi^*\sigma} (V),\qquad K^1_{\sigma+W_3(V)}(M)
 \cong K^1_{\pi^*\sigma} (V).
$$
Suppose that  $\pi: V\to M$ is an oriented real
vector bundle over $M$ of odd rank equipped with a positive definite
quadratic form.  Then there is a natural isomorphism
$$
K^0_{\sigma+W_3(V)}(M) \cong K^1_{\pi^*\sigma} (V), \qquad K^1_{\sigma+W_3(V)}(M)
 \cong K^0_{\pi^*\sigma} (V).
$$
\end{theorem}

We conclude by checking that we get the expected result for the Atiyah-Singer index theory for 
 Dirac operators on  Clifford modules and we explain how our result
  is relevant to string theoretic considerations.

 \section{Review of twisted K-theory}

We use \cite{AS} as our primary reference for twisted K-theory. We now 
review  the definition and some basic properties.

As noticed in \cite{AS}, on the one hand, for a smooth fiber bundle $f: X\to Y$ as appearing in
the push-forward map in the Introduction, the Hilbert bundle over $Y$ whose
fiber is the space $\cH$ of $L^2$-half densities along the fiber of $f$ does not
admit $U(\cH)$ with norm topology as structure group.
On the other hand for a compact Lie group, as in the definition of equivariant
twisted K-theory,  the regular representation
$G \to U(L^2(G))$ is not norm-continuous,   so for a stable $G$-Hilbert space
$\cH_G \cong \cH_G \otimes L^2(G)$, the action of $G$ on the space of Fredholm operators
is not norm-continuous. For these two reasons, we use  the compact-open topology
developed in the Appendix 1 in \cite{AS} on the unitary group and the projective
unitary group of a Hilbert space.

 Fix a standard $\ZZ_2$ graded Hilbert space
 $$\hat{\cH} =\hat{\cH}^+ \oplus\hat{\cH}^- $$
 such that
both $\hat{\cH}^+$ and $\hat{\cH}^-$ are infinite dimensional.  Let
$ Fred^{(0)}(\hat{\cH})$ be the space of self-adjoint Fredholm operators
$\tilde{A}$
in $\hat{\cH}$ of degree 1 with respect to the grading
such that $\tilde{A}^2$ differs from the identity by a compact operator.
The topology comes from the embedding
$\tilde{A} \mapsto (\tilde{A}, \tilde{A}^2-1)$
in $\cB(\hat{\cH}) \times \cK(\hat{\cH})$, where $\cB(\hat{\cH})$ is the space of bounded operators in $\hat{\cH}$
with the compact-open topology
and $\cK(\hat{\cH})$ is the Banach algebra of compact operators with norm
topology. Then as shown in \cite{AS},  $Fred^{(0)}(\hat{\cH})$ is a representing space
for $K$-theory, and admits a continuous conjugation action of   the projective unitary group
$PU(\hat{\cH})$ equipped with the compact-open
topology.   Note also that $PU(\hat{\cH})$ with the compact-open
topology acts continuously by conjugation on the Banach space  $\cK(\hat{\cH})$
of compact operators.

Given a twisting
\[
\sigma \in H^1(M, \ZZ_2) \oplus H^3(M, \ZZ),
\]
there is a unique isomorphism class of locally trivial bundles  over $M$ with fibre
an infinite
dimensional, separable, complex projective Hilbert space
in which a unitary involution is given in each fibre
\[
\hat{\cH} =\hat{\cH}^+ \oplus\hat{\cH}^-. 
\]
The structure group is 
the graded projective unitary group $PU_{gr}(\hat{\cH})$ with the compact-open topology,
where $PU_{gr}(\hat{\cH}) = U_{gr}(\hat{\cH})/U(1)$ and $U_{gr}(\hat{\cH})$ has
two connected components with the identity component 
\[
U^0_{gr}(\hat{\cH}) =  U(\hat{\cH}^+) \oplus  U(\hat{\cH}^-). 
\]
It was proved in \cite{Par1}\cite{Par2} that the grading preserving 
isomorphism classes of graded projective Hilbert bundles are classified by 
\[
H^1(M, \underline{PU_{gr}(\hat{\cH}) }) \cong H^1(M, \ZZ_2) \oplus H^3(M, \ZZ).
\]

Given   a twisting
$\sigma \in H^3(M, \ZZ)$, there is an associated isomorphism
class of projective
Hilbert space bundles (see \cite{AS} and \cite{BCMMS}).
The natural way to get a projective Hilbert space bundle with involution is
by introducing the $\ZZ_2$ graded space
$$\hat{\cH}= \cH \oplus \cH = \cH \otimes \CC^2,$$
 such that the structure group reduces to $U^0_{gr}(\hat{H})/U(1)$, hence
 the corresponding twisting has trivial components in  $H^1(M, \ZZ_2)$. For simplicity,
 throughout this paper, we only consider twistings   in $H^3(M, \ZZ)$. 

Fix a choice of such a projective Hilbert bundle $P_\sigma$
with involution associated to
$\sigma \in H^3(M, \ZZ)$ such that
there exists a good cover $\{U_\alpha\}$ of $M$ satisfying
\ba\label{trivial}
 P_\sigma |_{U_\alpha} \cong U_\alpha \times  \PP ( \cH)
\na
and such that the transition between two such trivializations is
given by a continuous function
\ba\label{transition}
g_{\alpha\beta}: U_\alpha \cap U_\beta  \longrightarrow PU(\cH),
\na
for $ PU(\cH)$ with either the norm topology or the compact-open topology,
Cf. Proposition 2.1 in \cite{AS}.

Denote by $\hat{P}_\sigma$ the graded tensor product
\[
\hat{P}_\sigma = P_\sigma \otimes \PP(\hat{\cH}).
\]
With respect to the local trivialization (\ref{trivial}), the transition functions are given by
$$\hat{g}_{\alpha\beta}: U_\alpha \cap U_\beta  \longrightarrow PU(\cH\otimes \hat{\cH}).$$
Denote by $\cP_\sigma$ the associated  principal bundle over $M$, and
 $Fred^{(0)}(\hat{P}_\sigma)$  the associated bundle of Fredholm operators
with fibre  $Fred^{(0)}(\cH\otimes \hat{\cH})$.

\begin{definition}\label{twisted:K}
 If $M$ is compact, then the twisted K-group
$K^0_{P_\sigma}(M)$ is defined to be the space of homotopy classes of sections
of $Fred^{(0)}(\hat{P}_\sigma)$, that is,
\[
K^0_{P_\sigma}(M) = [ M,  Fred^{(0)}(\hat{P}_\sigma)] \cong
[\hat{\cP}_\sigma, Fred^{(0)}(\hat{\cH})]_{PU(\hat{\cH})},
\]
where $[\hat{\cP}_\sigma, Fred^{(0)}(\hat{\cH})]_{PU(\hat{\cH})}$ is the group of  homotopy classes
of equivariant maps.
  If $M$ is locally compact, then the twisted K-group
$K^0_{P_\sigma}(M)$ is defined to be the space of homotopy classes of admissible sections
of $Fred^{(0)}(\hat{P}_\sigma)$, where a section $s$ is called admissible if there
is a compact set $K \subset M$ such that $s|_{M-K}$ is an invertible section of
$K^0_{P_\sigma}(M-K)$.
If $M_1$ is a closed subset of $M$, we define the relative twisted K-group
$$K^0_{P_\sigma}(M, M_0)= K^0_{P_\sigma}(M-M_0).$$
\end{definition}

In terms of  local trivializations (\ref{trivial})
and the corresponding transition functions, an element in
$K^0_{P_\sigma}(M)$ can be represented by a twisted family of local
sections:
\[
s_\alpha: U_\alpha \longrightarrow Fred^{(0)}(\hat{\cH}),
\]
such that $s_\beta = g_{\alpha\beta} s_\alpha  g_{\alpha\beta}^{-1}.$ From the definition,
we know that if $\sigma$ is trivial, then $K^0_{P_\sigma}(M)$ is isomorphic to
the ordinary topological K-group $K^0(M)$.

As discussed in \cite{AS}, there is an additive structure in $K^0_{P_\sigma}(M)$
given by the  operation of fiberwise direct sum after a choice  of an
isomorphism between $\hat{\cH}\otimes \CC^2$ and $\hat{\cH}$
(the space of these isomorphisms
is a contractible, hence, this is well-defined).
The group $K^0_{P_\sigma}(M)$ is functorial with respect to the pair $(M, P)$ and an
isomorphism of projective Hilbert bundles induces an isomorphism on their
twisted K-cohomology groups. In particular, $K^0_{P_\sigma}(M)$ admits a
natural action of $Aut(P_\sigma) \cong H^2(M, \ZZ)$.
 Up to a canonical isomorphism, the twisted K-group  $K^0_{P_\sigma}(M)$ depends only
 on the Dixmier-Douady class $\sigma$. With these understood, we
sometimes denote  the twisted K-group  $K^0_{P_\sigma}(M)$ by  $K^0_{\sigma}(M)$.

 In \cite{AS}, twisted groups $K^n_{P_\sigma}(M)$ for all $n\in \ZZ$
 are also defined. The bundle $Fred^{(0)}(\hat{P}_\sigma)$ has a
 base point on each fiber, represented by a chosen
 fiberwise  identification
 $(\hat{P}_\sigma)_x^+ \cong (\hat{P}_\sigma)_x^-$. Thus the fiberwise
 iterated based loop space forms a bundle
$$\Omega_X^nFred^{(0)}(\hat{P}_\sigma), $$
 whose fiber at $x$ is $\Omega^nFred^{(0)}((\hat{P}_\sigma)_x)$.
 Let the homotopy classes of sections of $\Omega_X^n Fred^{(0)}(\hat{P}_\sigma)$ be denoted
 $K^{-n}_{P_\sigma}(M)$. As shown in \cite{AS}, there exists a fiberwise homotopy equivalence
 $$Fred^{(0)}(\hat{\cH}) \longrightarrow \Omega^2 Fred^{(0)}(\hat{\cH})$$
 which is $PU(\hat{\cH})$-equivariant. Therefore, these groups
 $\{K^{-n}_{P_\sigma}(M)\}_{n\le 0}$ are periodic in $n$ with period $2$. This periodicity
 can be used to define twisted groups $K^n_{P_\sigma}(M)$ for all $n\in \ZZ$
  such that the twisted K-theory
 forms a  periodic cohomology theory
 of period 2 on the category of manifolds equipped with a
projective Hilbert bundle.
 In particular, an element in $K^1_{P_\sigma}$ can be
 represented by a   based $S^1$-family of local twisted sections
 \[
 \tilde{s}_\alpha : U_\alpha \times S^1 \longrightarrow  Fred^{(0)}(\hat{\cH})
 \]
 such that $\tilde{s}_\beta  = g_{\alpha\beta}\tilde{s}_\alpha  g_{\alpha\beta}^{-1}.$

  For a  manifold $M \times \RR$ with a twisting
given by the pull-back of a twisting $\sigma$
 on $M$, then
 \ba\label{0:1}
 K^0_{P_\sigma} (M\times \RR) \cong K^1_{P_\sigma} (M),
 \na
 and
 \ba\label{1:0}
 K^1_{P_\sigma} (M\times \RR) \cong K^{-1}_{P_\sigma} (M\times \RR)
 \cong K^{-2}_{P_\sigma} (M) \cong K^0_{P_\sigma} (M).
 \na

 Twisted K-theory satisfies the following basic properties  \cite{AS}.
 \begin{enumerate}
 \item Given two projective Hilbert bundles  $P_1$ and $P_2$ with involution  over $M$
 with their twisting  in
 $$ H^1(M, \ZZ_2) \oplus H^3(M, \ZZ),$$
 denote by $P_1\otimes P_2$ the graded tensor product of  $P_1$ and $P_2$.
 There is  a cup product  homomorphism
\ba\label{twisted:cup-product}
K^i_{P_1}(M)\times K^j_{P_{2}}(M) \longrightarrow
K^{i+j}_{P_1\otimes P_2}(M), 
\na
coming from the map $(A, A') \mapsto A\otimes 1 + 1\otimes  A'$ (the graded tensor product)
defined on the space of self-adjoint degree 1 Fredholm operators. 
  In particular,
 $K^0_{P_\sigma}(M)$ is a module
 over the untwisted K-group $K^0(M)$, extending the action of
 $H^2(M, \ZZ) \subset K^0(M)$.
 \item For any proper continuous map $f: M\to N$
  there exists a natural pull-back homomorphism,  which is a homotopy invariant of $f$,
 \[
 f^*: K^i_{f^*P_\sigma}(N) \longrightarrow  K^i_{P_\sigma}(M),
 \]
 where $f^*P_\sigma$ is the pull-back bundle of $P_\sigma$. In particular, if $M $
 is a closed subset of $N$ with the inclusion map $\iota: M\to N$, the
 pull-back map is just the restriction map.
 \item If $M$ is covered by two closed subsets $M_1$ and $M_2$ whose interiors cover $M$,
 for a projective Hilbert bundle $P$ over $M$,
 there is a Mayer-Vietoris  exact sequence
 \ba\label{Mayer-V1}
 \xymatrix{
 K^0_{P} (M)\ar[r]& K^0_{P_1} (M_1)\oplus  K^0_{P_2} (M_2)\ar[r]& K^0_{P_{12}} (M_1\cap M_2)
 \ar[d]\\
  K^1_{P_{12}} (M_1\cap M_2)\ar[u] & K^1_{P_1} (M_1)\oplus  K^1_{P_2} (M_2)\ar[l]&
   K^1_{P} (M) \ar[l]
   }
   \na
   where $P_1$, $P_2$ and $P_{12}$ are the restrictions of $P$ to
      $M_1$, $M_2$ and $M_1\cap M_2$ respectively.
\item If $M$ is covered by two  open subsets $U_1$ and $U_2$,
 for a projective Hilbert bundle $P$ over $M$,
  there is a Mayer-Vietoris  exact sequence
 \ba\label{Mayer-V2}
 \xymatrix{
 K^0_{P} (M)\ar[r]& K^1_{P_{12}} (U_1\cap U_2)\ar[r]& K^1_{P_1} (U_1)\oplus K^1_{P_2}( U_2)
 \ar[d]\\
 K^0_{P_1} (U_1) \oplus K^0_{P_2}( U_2)\ar[u] & K^0_{P_{12}} (U_1\cap U_2) \ar[l]&
   K^1_{P} (M) \ar[l]
   }
   \na
   where $P_1$, $P_2$ and $P_{12}$ are the restrictions of $P$ to
   $U_1$, $U_2$ and $U_1\cap U_2$ respectively.
      \end{enumerate}

\begin{remark} \begin{enumerate}
\item   The torsion subgroup of $H^3(M, \ZZ)$ can be
identified with the unitary Brauer group of $M$, the isomorphism classes of
finite dimensional projective bundles over $M$. For each
  torsion class  $\sigma$ in $H^3(M, \ZZ)$, K-theory with local coefficients is defined
  in \cite{DK} to get a complete Thom isomorphism in ordinary K-theory.
   We shall see that this is a prototype for the Thom isomorphism in Section 3.
\item (Cf. \cite{AS} \cite{Ros})
 Using the continuous action of $PU( \hat{\cH})$ with the compact-open topology
on $\cK(\hat{\cH})$, we  associate to $P_\sigma$ the bundle  $\cK_{P_\sigma}$ of
non-unital algebras whose fiber is the compact operators $\cK(\hat{\cH})$. Then
 the  twisted K-group
$K^0_{P_\sigma}(M)$ is canonically isomorphic to the K-theory of the
Banach algebra of sections of  $\cK_{P_\sigma}$, if $M$ is only locally compact, sections
of $\cK_{P_\sigma}$ are required to have compact support.
\item For an open embedding $\iota: U  \subset M$, the natural extension defines
a homomorphism
\[
K^i_{P|_U} (U)  \longrightarrow  K^i_{P} (M),
\]
which is  the push-forward map $\iota_*$ defined in next Section. Moreover, there is
the following exact sequence
\[
\xymatrix{
 K^0_{P|_U} (U)\ar[r]& K^0_{P} (M)\ar[r]& K^0_{P|_{M-U}} (M-U)
 \ar[d]\\
  K^1_{P|_{M-U}} (M-U)\ar[u] & K^1_{P } (M )  \ar[l]&    K^1_{P|_U} (U) \ar[l],
   }
 \]
 from which the Mayer-Vietoris exact sequences  (\ref{Mayer-V1})
 and  (\ref{Mayer-V2}) can  deduced for two pairs of  open
 embeddings $M_2-(M_1\cap M_2) \subset M_2$ and $M_2-(M_1\cap M_2)\subset M$ for
 two closed covering subsets $M_1$ and $M_2$.
\end{enumerate}
\end{remark}

Throughout this paper, we only deal with twisted K-theory over a locally compact
space.

\section{Thom isomorphism in twisted K-theory}\label{Thom}

\subsection{Twisted K-theory in the torsion case}\label{sec:torsion}

Given a torsion class $\sigma \in H^3(M, \ZZ)$, denote by
$\cP_\sigma$ the corresponding principal $PU(n)$-bundle over
$M$ with Dixmier-Douady invariant equal to $\sigma$. Let
$\Gamma_\sigma$ be the lifting bundle gerbe \cite{Mur}
\ba\label{bg:sigma}
\xymatrix{
 \Gamma_\sigma \ar[d]& \\
 \cP_\sigma  ^{[2]}\ar@< 2pt>[r]^{\pi_1} \ar@< -2pt>[r]_{\pi_2}
  & \cP_\sigma    \ar[d]\\
   &M}
\na
 over $M$ associated to $\cP_\sigma$ and the central extension
\[
1\to U(1) \longrightarrow U(n) \longrightarrow PU(n) \to 1.
\]
Note that $\Gamma_\sigma$ is the natural groupoid $U(1)$-extension of the groupoid
$\cP_\sigma^{[2]} = \cP_\sigma \times_M \cP_\sigma$ with the
source map given by $\pi_1: (y_1, y_2) \mapsto y_2$ and the target
map given by $\pi_2: (y_1, y_2)\mapsto y_1$.

\begin{definition} A bundle gerbe module $E$ of $\Gamma_\sigma$
(defined in \cite{BCMMS}), also called  a $\Gamma_\sigma$-module
in \cite{MurSin},
 is really a complex vector bundle $E$ over $\cP_\sigma$ with
a groupoid action of $\Gamma_\sigma$, i.e., an isomorphism
\[
\phi: \Gamma_\sigma \times_{(s,\pi)} E \longrightarrow E
\]
where $\Gamma_\sigma \times_{(s,\pi)} E$ is the fiber product of
the source $s: \Gamma_\sigma \to \cP_\sigma$  and $\pi: E\to
\cP_\sigma$ such that
\begin{enumerate}
\item $\pi\circ \phi (g, v) = t(g)$ for $(g, v) \in  \Gamma_\sigma
\times_{(s,\pi)} E$, and $t$ is the target map of $\Gamma_\sigma$.
\item $\phi$ is compatible with the bundle gerbe multiplication
$m: \Gamma_\sigma \times_{(s,t)}\Gamma_\sigma \to \Gamma_\sigma$,
which means
\[
\phi \circ (id \times \phi) = \phi\circ (m\times id).
\]
\end{enumerate}
\end{definition}

 Using the source map $r: \Gamma_\sigma \to \cP_\sigma$,
it is easy to see that $\Gamma_\sigma$ is a topologically trivial
$U(n)$-bundle over $\cP_\sigma$. If we think of $U(n)$ as acting on
$\cP_\sigma$ via the $PU(n)$-action (and using the quotient map $U(n)\to PU(n)$),
then
 we know that \cite{BCMMS} any $U(n)$-module $V$ gives rise to a bundle
 gerbe module $$\underline{V} = \cP_\sigma\times V$$ with $U(n)$-action
 $g\cdot (p, v) = (pg^{-1}, gv)$.

 The weight  of a $U(n)$-equivariant vector bundle is defined to be the 
 weight of the induced $U(1)$ action. Then
the additive  category of isomorphism classes of $\Gamma_\sigma$-modules $M$ is
equivalent to the additive  category of $U(n)$-equivariant vector
bundles over $\cP_\sigma$ of weight $1 \ mod (k)$, where  $k$ is
the order of $\sigma$. Therefore,
we conclude  that the K-group of bundle gerbe modules $K_{bg}(M,
\Gamma_\sigma)$, defined to be the Grothendieck group of the
additive category of $\Gamma_\sigma$-modules as in \cite{BCMMS}, is
isomorphic to the submodule $K^0_{U(n), (1)}(\cP_\sigma)$ of
weight $1 \ mod (k)$ in the  equivariant K-group
$K^0_{U(n)}(\cP_\sigma)$ regarded as an $R(U(1))$-module, where $R(U(1)$ denotes 
the complex representation ring of $U(1)$.

Here we recall that
$K^0_{U(n)}(\cP_\sigma)$ is the Grothendieck group of the additive
category of $U(n)$-vector bundles over the $U(n)$-manifold
$\cP_\sigma$.  The tensor product of $U(n)$-vector bundles  induces
a commutative ring structure in $K^0_{U(n)}(\cP_\sigma)$. The
natural morphism from the representation ring  $R(U(n))$ to
$K^0_{U(n)}(\cP_\sigma)$ defines a $R(U(n))$-module structure on
$K^0_{U(n)}(\cP_\sigma)$. The $R(U(1))$-module structure on
$K^0_{U(n)}(\cP_\sigma)$ is defined by the action of $U(1)$ on
$U(n)$-vector bundles.  

The following proposition gives a
correct  proof of  Proposition 6.4 of \cite{BCMMS}.
The discussion in \cite{BCMMS} does not handle
$U(n)$-equivariant contractibility of the unitary group of Hilbert space correctly
(and the argument below fixes this using the appendix to \cite{AS}), and the index map 
\[
ind:  [ \cP_\sigma, Fred (C^n\otimes \cH)]_{U(n)}\longrightarrow K_{bg}(M, \Gamma_\sigma),
\]
constructed in section 6.2 of \cite{BCMMS} is neither injective nor surjective, contrary to the
claim in  the proof of Proposition 6.4 of \cite{BCMMS}.

\begin{proposition}\label{torsion:twist=bg}
For a torsion  class $\sigma \in H^3(M, \ZZ)$ of order $k$, let $\cP_\sigma$ be
a principal $PU(n)$-bundle whose Dixmier-Douady class is $\sigma$, then
the twisted K-group $K^0_{\cP_\sigma}(M)$ is isomorphic to   $K^0_{U(n), (1)}(\cP_\sigma)$, hence
\[
K^0_{\cP_\sigma}(M) = K_{bg}(M, \Gamma_\sigma).
\]
\end{proposition}
\begin{proof}
Using the representing space of $K_{U(n)}^{0}$ developed in Appendix 3 of \cite{AS}, the proof
is rather straight forward. Let $\cH_{U(n)}$ be a stable $U(n)$-Hilbert space, then
the closed subspace of $U(n)$-continuous Fredholm operators in $Fred(\cH_{U(n)})$
\[
Fred_{cts}(\cH_{U(n)})= \{D \in Fred(\cH_{U(n)}): g\mapsto gDg^{-1} \text{  is continuous  for}
 \  g\in U(n) \},
\]
with norm topology, is a representing space of $K_{U(n)}^{0}$. That means,
\[
K^0_{U(n)}(\cP_\sigma) = [ \cP_\sigma, Fred_{cts}(\cH_{U(n)})]_{U(n)},
\]
the group of all homotopy classes of equivariant maps. As a $R(U(1))$-module, the
submodule of  weight $1 \ mod (k)$
\[K^0_{U(n), (1)}(\cP_\sigma)=
[ \cP_\sigma, Fred_{cts}(\cH_{U(n), (1)})]_{U(n)},
\]
where $Fred_{cts}(\cH_{U(n), (1)})$ is the $U(n)$-Hilbert subspace of $\cH_{U(n)}$
of weight $1 \ mod (k)$ under the action of $U(1)$.  By the results in Appendix 3 in \cite{AS},
we obtain that
\[
K^0_{U(n), (1)}(\cP_\sigma)= [ \cP_\sigma, Fred^{(0)}(\cH_{U(n),
(1)}) ]_{U(n)}.\] From the definition of the twisted K-group,
applying the action of $U(n)$ on $\CC^n\otimes \hat{\cH}$ defined
by $g\otimes 1$, we know that
\[
\begin{array}{lll}
K^0_{\cP_\sigma}(M)
   &=&
 [ \cP_\sigma\times_{PU(n)}PU(\CC^n\otimes \hat{\cH}),
  Fred^{(0)}(\CC^n\otimes \hat{\cH})]_{PU(\CC^n\otimes \hat{\cH})}\\[2mm]
  &\cong &   [ \cP_\sigma\times_{PU(n)}PU(\cH_{U(n), (1)}),
    Fred^{(0)}(\cH_{U(n), (1)})]_{PU(\cH_{U(n), (1)})}\\[2mm]
  &\cong &  [ \cP_\sigma, Fred^{(0)}(\cH_{U(n), (1)}) ]_{PU(n)}\\[2mm]
  &\cong & K^0_{U(n), (1)}(\cP_\sigma)\\[2mm]
  &\cong &K_{bg}(M, \Gamma_\sigma) .
 \end{array} \]
  \end{proof}

 Let $\cE_{bg}(M, \Gamma_\sigma)$ be the additive category of
 $\Gamma_\sigma$-modules.  Let $\cA_\sigma$ be the 
  Azumaya bundle over $M$ with Dixmier-Douady class $\sigma$
  \[
  \cA_\sigma = \cP_\sigma \times_{Ad} M_n(\CC),
  \]
  and where $Ad$ is the adjoint representation of $PU(n)$ on $M_n(\CC)$
  by conjugation.  An $\cA_\sigma$-module is a complex vector
  bundle over $M$ with a fiberwise $\cA_\sigma$ action
   $$
   \cA_\sigma \times_M E \longrightarrow E.
   $$
   Applying the  natural representation of $U(n)$ on $\CC^n$, then
 $\underline{\CC^n}^*= \CC^n \times \cP_\sigma$ is a  $\Gamma_{-\sigma}$-module.
 Tensoring with $\underline{\CC^n}^*$, there is
  a natural equivalence of categories between $\cE_{bg}(M, \Gamma_\sigma)$
  and $\cE^{\cA_\sigma}(M)$, the category of $\cA_\sigma$-modules.
  Therefore, we obtain
  \ba\label{twisted:azumaya}
  K^0_{\cP_\sigma}(M) \cong K (\cE^{\cA_\sigma}(M)),
  \na
  the Grothendieck group of the category of $\cA_\sigma$-modules.

\subsection{Thom isomorphism in twisted K-theory for the torsion case}

Consider an  oriented real vector bundle $V$ of even rank over $M$
with a positive definite quadratic form $Q_V$. Denote by $\cFr$
the frame bundle of $V$, the principal $SO(2n)$-bundle of
positively oriented orthonormal frames, i.e.,
$$
V= \cFr\times_{\rho_{2n}} \RR^{2n}, $$
where $\rho_n$ is the standard representation of $SO(2n)$ on $\RR^n$.
 The lifting bundle gerbe (Cf. \cite{Mur}) associated to the frame bundle and
the central   extension
  \[
  1\to U(1) \longrightarrow Spin^c(2n) \longrightarrow SO(2n) \to 1
  \]
  has its Dixmier-Douady class given by  the integral third Stiefel-Whitney class
   $W_3(V)\in H^3(M, \ZZ)$.  We denote this lifting bundle
  gerbe by $\Gamma_{W_3(V)}$
\[
  \xymatrix{
 \Gamma _{W_3(V)}\ar[d]& \\
 \cFr  ^{[2]}\ar@< 2pt>[r]^{\pi_1} \ar@< -2pt>[r]_{\pi_2}
  & \cFr    \ar[d]\\
   &M.}
\]
Note that $\Gamma_{W_3(V)} \cong\Gamma_{spin} \times_{ \ZZ_2} U(1) $ where
$\Gamma_{spin}$ is a $\ZZ_2$-bundle gerbe associate to $\cFr \to M$ and the central
extension
\[
  1\to \ZZ_2 \longrightarrow Spin (2n) \longrightarrow SO(2n) \to 1.
  \]
The Clifford bundle  $Cl(V)$ of $V$ (\cite{LM}) is defined to be
\[
Cl(V) = \cFr\times_{cl(\rho_{2n})}Cl(\RR^{2n}),
\]
where $Cl(\RR^{2n})$ is the Clifford algebra of $\RR^n$ with the standard positive
definite quadratic form, and $cl(\rho_{2n})$ is the induced representation of  $SO(2n)$
on $Cl(\RR^{2n})$. The fiberwise Clifford multiplication in $Cl(V)$, as a continuous map,
\[
Cl(V)\times_M Cl(V) \longrightarrow Cl(V)
\]
defines an algebra structure on $Cl(V)$. Hence $Cl(V)$ is called the
bundle of Clifford algebras. A complex vector bundle $E$ over $M$ is called a complex
$C(V)$-module if there is a continuous map, called the Clifford action,
 $$
 Cl(V) \times_M E \longrightarrow  E
 $$ such that
each fiber $E_x$ ($x\in M$) is a  complex  $Cl(V_x)$-module.
Denote  by $Cl_\CC(V)= Cl(V)\otimes \CC $  the complexified Clifford bundle.  Note that
a complex $C(V)$-module is a $Cl_\CC(V)$-module.

  Let $\cE^V(M)$ be the category of complex  $C(V)$-modules whose morphisms are
  vector bundle morphisms which commute with   the $Cl(V_x)$-actions over each point $x\in M$.
  Let $\cE (M, \Gamma_{W_3(V)})$ be the category of $\Gamma_{W_3(V)}$-modules
  whose morphisms are vector bundle morphisms which commute with the
  groupoid action of $\Gamma_{W_3(V)}$. Under the Whitney sum
  of vector bundles, $\cE^V(M)$  and $\cE (M, \Gamma_{W_3(V)})$
  are additive categories. As observed in \cite{MurSin}, we now show that
  these two categories   are  in fact equivalent.

  \begin{lemma}\label{Cliff=gb} There is a natural functor
  \[
 \Psi:  \cE^V(M) \longrightarrow \cE (M, \Gamma_{W_3(V)}),
  \]
   which is a category equivalence. In particular, the Grothendieck
   group of $\Gamma_{W_3(V)}$-modules is isomorphic to
the Grothendieck group of $\cE^V(M)$ which is denoted by
$K(\cE^V(M))$
\[
K_{bg}(M, \Gamma_{W_3(V)}) \cong K(\cE^V(M)).
\]
 \end{lemma}
 \begin{proof}
 There exists a distinguished $\Gamma_{W_3(V)}$-module given by restricting
  the irreducible   complex representation   of $Cl(\RR^{2n})$
  to $Spin (2n) \subset Cl(\RR^{2n})$. Denote   the unique  irreducible
   complex representation of $Cl(\RR^{2n})$ by
   \[
   \Delta_{2n}^\CC: Spin (2n)  \to GL_\CC (S_{2n}).
   \]
    Applying the groupoid structure of $\Gamma_{spin}$, the representation $\Delta_{2n}^\CC$
    defines  a natural $\Gamma_{spin}$-module structure on
    \[
    \underline{S_{2n}} = S_{2n} \times \cFr,
    \]
    see the discussion before Proposition \ref{torsion:twist=bg}. As $Spin^c(2n)
    = Spin(2n) \times_{\ZZ_2} U(1)$,  $S_{2n}$ is
    a $Spin^c(2n)$-module with $U(1)$-action given by scalar multiplication of
     $U(1)  \subset \CC$.  This
    implies that $\underline{S_{2n}}$ is  a $\Gamma_{W_3(V)}$-module. Note that the dual bundle
    $\underline{S_{2n}^*}$  is also a $\Gamma_{W_3(V)}$-module.

      Let $\pi: \cFr\to M$ be the frame bundle associated to $(V, Q_V)$.
    Then $Cl(\pi^* V)$ is  the Clifford bundle over $\cFr$ associated
    to $(\pi^*V, \pi^*Q_V)$.
    Given a  complex  $C(V)$-module $E$,  $\underline{S_{2n}}^*\otimes_{Cl(\pi^* V)} \pi^*E$
    is a vector bundle over $\cFr$ which
    admits a $\Gamma_{W_3(V)}$-module structure coming from $Spin^c(2n)$-action
   on $\underline{S_{2n}}$.
   This defines a natural functor
   \[
   \begin{array}{cccl}
   \Psi: &\cE^V(M) &\longrightarrow & \cE (M, \Gamma_{W_3(V)})\\[2mm]

  & E &\mapsto & W_E  =  \underline{S_{2n}}^* \otimes_{Cl(\pi^* V)}  \pi^*E,
   \end{array}
   \]
   such that $Hom_{Cl(V)}(E, F) \cong Hom_{\Gamma_{W_3(V)}}( W_E, W_F)$
   as $Cl_\CC (\pi^*V) \cong End (\underline{S_{2n}})$.

   Now we  show that  $\Psi$ is a functor setting up a category equivalence. We only need to construct
   an inverse to $\Psi$, up to isomorphism.  Given a $\Gamma_{W_3(V)}$-module
   $W$,  $W$ is a $Spin(2n)$-equivariant vector bundle over $\cFr$. Then
   $W\otimes \underline{S_{2n}}$ is an $SO(2n)$-equivariant vector bundle over $\cFr$,
   hence, descends to  a complex vector bundle over $M$, denoted by $E_W$. The Clifford action of
   $Cl(\RR^{2n})$-action on $\underline{S_{2n}}$ defines a complex $Cl(V)$-module
   structure on $E_W$.  Using the isomorphism $Cl_\CC (\pi^*V) \cong End (\underline{S_{2n}})$ again,
   we see that $\Psi (E_W) = W$. Morphisms in
    $\cE_{bg}(M, \Gamma_{W_3(V)})$ and $ \cE^{V}(M)$
    can also be identified.
   \end{proof}

   This identification of $\cE^V(M)$ with $\cE (M, \Gamma_{W_3(V)})$ is very important, as we can
   obtain some nice properties of $\cE^V(M)$ from $\cE (M, \Gamma_{W_3(V)})$. Using the
   Peter-Weyl theorem as in the proof of Proposition 2.4 in \cite{Seg}, for any
   $\Gamma_{W_3(V)}$-module $W$, there is a $\Gamma_{W_3(V)}$-module $W^\bot$ such
   that
   \ba\label{pseudo-abel}
   W \oplus W^\bot \cong \underline{ \CC^N} = \CC^N \times \cFr
   \na
   for a $Spin(2n)$-module $\CC^N$.  With this property, we can identify
   $K(\cE^V(M))$ with the Grothendieck group, denoted by $K^V(X)$, of
   the following  forgetful functor (see Ch IV.5 in \cite{Kar0})
   \[
   \cE^{V\oplus \underline{\RR}}(M)
   \longrightarrow  \cE^V(M)
   \]
   where $\underline{\RR}= M\times \RR$ is the trivial bundle of rank one with the
   standard quadratic form.   Let $B(V)$ (resp. $S(V)$ ) be the ball bundle (resp. the sphere
   bundle) of $V$.  Define
   \[
   K_n^V(M) = K^V (M\times B^n, M\times S^{n-1}).
   \]
   The following Thom isomorphism theorem  is established in \cite{Kar0}.

   \begin{theorem} (Theorems IV.5.11 and IV.6.21 \cite{Kar0}) There is a natural
   isomorphism
   $$K_n^V(M) \cong K_n (B(V), S(V)) = K_n(V)$$
   which generalizes the standard Thom isomorphism (i.e., tensoring with the
   Thom class of $V$) for $V$ admitting
   a $Spin^c$-structure (i.e., $W_3(V) =0$).
   \end{theorem}

   Applying Theorem III.4.12 in \cite{Kar0} to the category $\cE^V(M)$ which is
   a pseudo-abelian category (see page 28 in \cite{Kar0}) by (\ref{pseudo-abel}),
   we know that $K(\cE^V(M))$ is isomorphic to the Grothendieck group of
   the following forgetful functor
   \[
   \bigl(\cE^V(M)\bigr)^{Cl(\underline{\RR}) } \longrightarrow \cE^V(M),
   \]
   where $\bigl(\cE^V(M)\bigr)^{Cl(\underline{\RR}) }$ is a category
   whose objects are pairs  $(E, \rho)$, where
   $E\in \cE^V(M)$ and $\rho$ is an $\RR$-linear bundle homomorphism
   \[
     \rho:  Cl(\underline{\RR}) \longrightarrow End_{Cl(V)} (E).
     \]
     As $\bigl(\cE^V(M)\bigr)^{Cl(\underline{\RR}) }\cong \cE^{V\oplus \underline{\RR}}(M)$,
     we obtain
     \ba\label{KV=KE}
     K(\cE^V(M) ) \cong K^V(M) \cong K (V).
     \na
     By substitution $M\mapsto (M\times B^n, M\times S^{n-1})$, we have
     \[
     K_n(\cE^V(M)) \cong   K_n (V).
     \]

Together with Proposition \ref{torsion:twist=bg}, we have obtained
the following isomorphisms:
\[
\begin{array}{lll}
K_{W_3(V)}^0(M) &\cong& K_{bg}^0(M, \Gamma_{W_3(V)})\\[2mm]
&\cong &  K(\cE^V(M))\\[2mm]
&\cong & K^V(M) \cong K(V).
\end{array}\]
This  is the  Thom isomorphism in K-theory for any real oriented
vector bundle established by Donovan-Karoubi in \cite{DK} using
K-theory with local coefficients.

Given a torsion element $\sigma \in H^3(M, \ZZ)$, let $\cP_\sigma$
be the associated principal $PU(n)$-bundle and  $\Gamma_\sigma$ be
the corresponding  lifting bundle gerbe given by (\ref{bg:sigma}).
 We now describe the Thom isomorphism in twisted K-theory for the 
torsion case.

\begin{theorem}\label{Thom0} Suppose that  $\pi: V\to M$ be an oriented real
vector bundle over $M$ of even rank with positive definite
quadratic form.  There is a natural isomorphism
$$
K^0_{\sigma+W_3(V)}(M) \cong K^0_{\pi^*\sigma} (V), \qquad K^1_{\sigma+W_3(V)}(M)
 \cong K^1_{\pi^*\sigma} (V).
$$
 Suppose that $\pi: V\to M$ be an oriented real
vector bundle over $M$ of odd rank with positive definite
quadratic form.  There is a natural isomorphism
$$
K^0_{\sigma+W_3(V)}(M) \cong K^1_{\pi^*\sigma} (V), \qquad K^1_{\sigma+W_3(V)}(M)
 \cong K^0_{\pi^*\sigma} (V).
$$
\end{theorem}
\begin{proof}  We know that   $K_{\sigma+W_3(V)}(M) \cong
 K(\cE(M, \Gamma_\sigma\otimes \Gamma_{W_3(V)} )$, the
 Grothendieck group of the additive category of $\Gamma_\sigma\otimes
 \Gamma_{W_3(V)}$-modules (complex vector bundles over
 $\cFr\times_{M}\cP_\sigma$ with the groupoid action of $\Gamma_\sigma\otimes
 \Gamma_{W_3(V)}$). 
   Then the category of $\Gamma_\sigma\otimes
 \Gamma_{W_3(V)}$-modules is isomorphic to $\bigl(\cE^V(M)\bigr)^{\cA_\sigma}$,  the category
of $Cl(V)$-modules admitting an $\cA_\sigma$-action. As $K(\cE^V(M))
\cong K(V)$, it is straight forward to see that
$K(\bigl(\cE^V(M)\bigr)^{\cA_\sigma}) \cong K_{\pi^*\sigma}(V)$
for, we just apply the explicit isomorphism between $K(\cE^V(M)) $
and $ K(V)$ (Cf. page 225 in \cite{Kar0}), in particular, if $E\in
\bigl(\cE^V(M)\bigr)^{\cA_\sigma}$, then $\pi^* (E) \in
\cE^{\pi^*\cA_\sigma}(V)$. This implies that
$$
K_{\sigma+W_3(V)}(M) \cong K_{\pi^*\sigma} (V).
$$
As $K^1_{\sigma+W_3(V)}(M)\cong K_{\sigma+W_3(V)}(M\times \RR)$,
and $K_{\sigma+W_3(V)}(M\times \RR)  \cong K^0_{\pi^*\sigma} (V \times\RR)$, we
obtain $$K^1_{\sigma+W_3(V)}(M)\cong K^1_{\pi^*\sigma} (V).$$

For  an oriented real vector bundle $V$ over $M$ of odd rank with positive definite
quadratic form, we apply the Thom isomorphism to $V\oplus \underline{\RR}$.
We have
\ba\label{iso:1}
K_{\sigma+W_3(V)}(M) \cong K_{\pi^*\sigma} (V\oplus \underline{\RR}),
\qquad K^1_{\sigma+W_3(V)}(M)
 \cong K^1_{\pi^*\sigma} (V\oplus \underline{\RR}).
 \na
Now we can think of $V\oplus \underline{\RR}$ as a rank one real vector bundle
over $V$, and  apply the Thom isomorphism again to get
\ba\label{iso:2}
K_{\pi^*\sigma} (V\oplus \underline{\RR}) \cong K^1_{\pi^*\sigma} (V),
\qquad K^1_{\pi^*\sigma} (V\oplus \underline{\RR}) \cong K_{\pi^*\sigma} (V).
\na
Put (\ref{iso:1}) and (\ref{iso:2}) together, we have
$$K_{\sigma+W_3(V)}(M) \cong K^1_{\pi^*\sigma} (V), \qquad K^1_{\sigma+W_3(V)}(M)
 \cong K^0_{\pi^*\sigma} (V).
$$
\end{proof}

\subsection{Thom isomorphism in twisted K-theory for the non-torsion case}\label{Thom:non-torsion}

 The virtue of our approach in the torsion case is that it may be adapted to the non-torsion
 case.
  For a non-torsion twisting $\sigma \in H^3(M, \ZZ)$, let $\Gamma_\sigma$
  be the lifting bundle gerbe associated to a $PU(\cH)$-principal bundle
  $$\pi_\sigma: \cP_\sigma \to M,$$ where $\cH$ is the standard $\ZZ_2$ graded separable complex
   Hilbert space. Note that $\Gamma_\sigma$ has a groupiod structure with
   the space of objects given by $\cP_\sigma$, it is a groupoid $U(1)$-extension
   of the natural groupoid $\cP_\sigma^{[2]} = \cP_\sigma \times_M \cP_\sigma$.

  Let $U_\cK$ be the  normal subgroup of $U(\cH)$ of unitary operators of the form
  $1 +$ compact, and $U_2$ be the normal subgroup of $U(\cH)$ of unitary operators of the form
  $1 +$ Hilbert-Schmidt operator. 
  
  \begin{definition} Denote by $\cE^{U_2}_{bg}(M, \Gamma_\sigma)$ the additive category of
  $\Gamma_\sigma$-modules with $U_2$ structure group. Here a $\Gamma_\sigma$-module $W$
   with $U_2$ structure is a Hilbert bundle  $W$ 
   over $\cP_\sigma$ with structure group $U_2$
   and an action of the groupoid $\Gamma_\sigma$,  i.e., an isomorphism
\[
\phi: \Gamma_\sigma \times_{(s,\pi)} W \longrightarrow W
\]
where $\Gamma_\sigma \times_{(s,\pi)} W$ is the fiber product of
the source $s: \Gamma_\sigma \to \cP_\sigma$  and $\pi: W\to
\cP_\sigma$ such that
\begin{enumerate}
\item $\pi\circ \phi (g, v) = t(g)$ for $(g, v) \in  \Gamma_\sigma
\times_{(s,\pi)} E$, and $t$ is the target map of $\Gamma_\sigma$.
\item $\phi$ is compatible with the bundle gerbe multiplication
$m: \Gamma_\sigma \times_{(s,t)}\Gamma_\sigma \to \Gamma_\sigma$,
which means
\[
\phi \circ (id \times \phi) = \phi\circ (m\times id).
\]
\end{enumerate}
\end{definition}

  In  \cite{BCMMS}, it was shown
  that 
  \[ K_\sigma (M) \cong [\cP_\sigma, Fred^{(0)}(\cH) ]^{PU(\cH)}
  \cong  K(\cE^{U_2}_{bg}(M, \Gamma_\sigma)).
  \]
   To see this, we need a $PU(\cH)$-equivariant homotopy
  equivalence between $Fred^{(0)}(\cH)$ and $\ZZ\times BU_\cK$ where
  $BU_\cK= U(\cH)/U_\cK$ with the conjugation action of $PU(\cH)$.  We
  apply  Bott periodicity:
  \[\begin{array}{rcl}
  U &\sim &\Omega (\ZZ \times BU ) \\[2mm]
   \ZZ \times BU & \sim &\Omega U
   \end{array}
   \]
   where $U=U_K,$ or $U_2$.
   Then use the  $PU(\cH)$-equivariant homotopy
  equivalence between $U_\cK$ and $U_2$,  to  obtain a
   $PU(\cH)$-equivariant homotopy
  equivalence between $BU_\cK$ and $BU_2$, where we choose  a homotopy model
  of $BU_2$  as the  component of  the  identity in   $\Omega U_2$,
  on which  $PU(\cH)$ acts by conjugation. This implies that
  \ba\label{twisted:BU_2}
  K_\sigma (M) \cong [\cP_\sigma, BU_2]^{PU(\cH)},
  \na
  that is, the twisted K-group of $K_\sigma (M)$ can be described by
  isomorphism classes of $\Gamma_\sigma$-modules with $U_2$ structure group (Cf.
  Proposition 7.2 in \cite{BCMMS})
  
  Recall that the unitary group $U(\cH)$ with the compact-open topology
  is equivariantly contractible and acts continuously by conjugation  on the Banach space
  $\cK(\cH)$ of compact operators in $\cH$, and also on the Hilbert space
  $$End_{HS}(\cH)= \cH\otimes \cH^*$$ of Hilbert-Schmidt operators (Cf. Appendix in \cite{AS}).

  \begin{definition} 
  Let $\cA_\sigma$ be the bundle of Hilbert-Schmidt operators  associated to $\cP_\sigma$,
      \[
  \cA_\sigma  = \cP_\sigma \times_{PU(\cH)}   End_{HS}(\cH),
  \]
  where $PU(\cH)$ acts on the space of Hilbert-Schmidt operators $End_{HS}(\cH)$ by conjugation.
   Let
   $\cE^{\cA_\sigma}_{U_2}(M)$ be the category of $\cA_\sigma$-modules with $U_2$ structure
   group (where
   a $\cA_\sigma$-module $E$ with $U_2$ structure group is a Hilbert bundle $E$
   over $M$ with structure group
   $U_2$ and an action
     \[
     \cA_\sigma\times_M E\longrightarrow E.
     \]
\end{definition}

This  bundle of Hilbert-Schmidt operators  plays the role that the 
 Azumaya bundle plays for the torsion case.  Just as
 there exists an equivalence functor between the category of
 bundle gerbe modules and the category of modules of  the 
 Azumaya bundle, we have the corresponding equivalence of categories
 between the category of
  $\Gamma_\sigma$-modules with $U_2$ structure group and the category of 
   $\cA_\sigma$-modules with $U_2$ structure
   group.

     \begin{lemma} \label{twsited=HS-module}
     There is a natural functor $\Psi:
     \cE^{U_2}_{bg}(M, \Gamma_\sigma) \to \cE^{\cA_\sigma}_{U_2}(M)$
     which defines an equivalence of categories. In particular,
     \[
     K_\sigma (M) \cong  K(\cE^{U_2}_{bg}(M, \Gamma_\sigma)) \cong K(\cE^{ \cA_\sigma}_{U_2}(M)).
     \]
     \end{lemma}
     \begin{proof}
     There exists a distinguished element $$\underline{\cH} = \cH \times \cP_\sigma$$
     in $\cE^{U_2}_{bg}(M, \Gamma_\sigma)$ and
     the dual bundle $\underline{\cH}^*$ is a
     $\Gamma_\sigma^*$-module    with $U_2$ structure group, where
     $ \Gamma_\sigma^*$ is the dual bundle gerbe associated to $\Gamma_\sigma$. Given a
     $\Gamma_\sigma$-module $W$  with $U_2$ structure group, then
     \[
     \underline{\cH}^* \otimes W
     \]
     is a Hilbert bundle over $\cP_\sigma\otimes \cP_{-\sigma} \cong M \times PU(\cH)$, whose
       structure group  has a reduction to $U_2$, see Section 4.2 in \cite{MatSte} for the similar
       argument for $U_\cK$.
     Notice that $\underline{\cH}^*\otimes W$ is a $PU(\cH)$-equivariant bundle
     over $\cP_\sigma$ and hence
     descends to a  Hilbert bundle over $M$ with structure group
   $U_2$.  Denote this bundle by $\Psi(W)$. Then $ \Psi(W)$ admits a natural
   action of $\cA_\sigma$, fiberwisely defined by the action of $End_{HS}(\cH)$ on $\cH$. Hence,
   $ \Psi(W)$ is a $\cA_\sigma$-module with $U_2$ structure group.

    Conversely, given a
   $\cA_\sigma$-module $E$ with $U_2$ structure group, then
   $$\underline{\cH} \otimes_{ \pi_\sigma^*\cA_\sigma}  \pi_\sigma^*E
   \cong Hom_{\pi_\sigma^*\cA_\sigma} (\underline{\cH}^*, \pi_\sigma^*E) $$
    is a
    $\Gamma_\sigma$-module with $U_2$ structure group satisfying
    $\Psi(\underline{\cH} \otimes_{ \pi_\sigma^*\cA_\sigma} \pi_\sigma^*E ) = E$. Here we use
    $$\pi_\sigma^*\cA_\sigma \cong \underline{\cH}^* \otimes \underline{\cH} ,$$
    and $\pi_\sigma^*E$ is a $\pi_\sigma^*\cA_\sigma$-module with $U_2$
    structure group.  Morphisms in
    $\cE^{U_2}_{bg}(M, \Gamma_\sigma)$ and $ \cE^{\cA_\sigma}_{U_2}(M)$
    can also be identified.  Therefore,
    $\Psi$ defines a category equivalence.
    \end{proof}

Now we can prove the Thom isomorphism in twisted K-theory for the general non-torsion
case.

\begin{theorem}\label{Thom1} Let $\sigma\in H^3(M, \ZZ)$ be a non-torsion element.
Suppose that $\pi: V\to M$ is an oriented real
vector bundle over $M$ of even rank with positive definite
quadratic form.  There is a natural isomorphism
$$
K^0_{\sigma+W_3(V)}(M) \cong K^0_{\pi^*\sigma} (V),\qquad K^1_{\sigma+W_3(V)}(M)
 \cong K^1_{\pi^*\sigma} (V).
$$
Suppose that  $\pi: V\to M$ is an oriented real
vector bundle over $M$ of odd rank with positive definite
quadratic form.  There is a natural isomorphism
$$
K^0_{\sigma+W_3(V)}(M) \cong K^1_{\pi^*\sigma} (V), \qquad K^1_{\sigma+W_3(V)}(M)
 \cong K^0_{\pi^*\sigma} (V).
$$
\end{theorem}
\begin{proof}
Let $\Gamma_\sigma$ and $\Gamma_{W_3(V)}$ be lifting bundle gerbes
associated to the corresponding $PU(\cH)$-bundles $\pi_\sigma: cP_\sigma \to M$ and
$\pi: \cFr\to M$. Then
\[\begin{array}{lll}
K_{\sigma+ W_3(V)}(M) &\cong & K(\cE^{U_2}_{bg}(M, \Gamma_\sigma\otimes\Gamma_{W_3(V)} ) )\\[2mm]
&\cong &  K(\cE^{ \cA_\sigma\otimes Cl(V)}_{U_2 }(M)) \\[2mm]
&\cong &
K(\cE^{\pi_\sigma^* \cA_\sigma }_{U_2}(V))\\[2mm]
&\cong & K_{\pi^*\sigma} (V).
\end{array}
\]
The Thom isomorphisms for other cases can be obtained by the same arguments as in the
proof of Theorem \ref{Thom0}.
\end{proof}

\begin{remark} For a locally compact and non-compact manifold $M$, some care is
 required as  twisted K-theory is defined
with compact support. We always assume that the isomorphism (\ref{twisted:BU_2})
 should be understood with respect to compactly supported  $PU(\cH)$-equivariant
 homotopy equivalent classes,   and
  $\Gamma_\sigma$-modules  in $\cE^{U_2}_{bg}(M, \Gamma_\sigma)$
  should be isomorphic to $\underline{\cH}$ away from compact sets.
  \end{remark}

\section{Push-forward map in twisted K-theory}
 
 In this section, we will establish the main theorem of this
 paper: the existence of the push-forward map in twisted K-theory for
 any differentiable  map $f: X \to Y$. We 
 equip $X$ and $Y$ with Riemannian metrics. So the Clifford bundles
 associated with their tangent bundles are well-defined.
 
  Denote by $W_3(f)$ the image of $w_2(X)- f^*(w_2(Y)) \in H^2(X, \ZZ_2)$ under
   the Bockstein homomorphism
  \[
  \beta: H^2(X, \ZZ_2) \longrightarrow H^3(X, \ZZ).
  \]
  For $f:X \to pt$, $W_3(f)$ is the third integer Stiefel-Whitney class $W_3(X)$ of
  $X$.
  
  For a twisting $\sigma \in H^3(Y, \ZZ)$, we choose a $PU(\cH)$-principal bundle
  $$\pi_\sigma: \cP_\sigma \to Y.$$  
  Let $\cA_\sigma$ be the bundle of Hilbert-Schmidt operators  associated to $\cP_\sigma$,
      \[
  \cA_\sigma  = \cP_\sigma \times_{PU(\cH)}   End_{HS}(\cH),
  \]
  where $PU(\cH)$ acts on the space of Hilbert-Schmidt operators $End_{HS}(\cH)$ by conjugation.
  
  The twisted K-group $K^*_{\sigma }(Y)$ is  defined
  with respect to this $PU(\cH)$-principal bundle  $ \cP_\sigma$ using Definition \ref{twisted:K},
  Equivalently, $K^*_{\sigma }(Y)$ can be defined as the Grothendieck group of the category of
  $\cA_\sigma$-modules with $U_2$-structure (see Lemma \ref{twsited=HS-module}). 
 
  Using the pull-back bundle $f^*\cA_\sigma$ tensored with the  Clifford bundle 
  $Cl( TX \oplus f^*TY)$ associated with
  $$TX \oplus f^*TY,$$ we get a bundle $\cA_{f^*\sigma + W_3(f)}$ 
  of  Hilbert-Schmidt operators  over $X$ with  the twisiting 
  $$f^*\sigma + W_3(f)    \in H^3(X, \ZZ), $$
   here the Hilbert space is the graded tensor product
   of $\cH$ with the unique $\ZZ_2$-graded irreducible module of
   the Clifford algebra. Then the  twisted K-group
  $K^*_{f^*\sigma + W_3(f)}(X)$ is defined to be the 
  Grothendieck group of the category of $\cA_{f^*\sigma+W_3(f)}$-modules with 
  $U_2$-structure.

\begin{theorem}\label{main:theorem}
Let  $f: X\to Y$  be a  differentiable map and $\sigma
  \in H^3(Y,\ZZ)$.
  There exists a natural push-forward map
  \[
  f_{!}: K^*_{f^*\sigma + W_3(f)}(X) \longrightarrow K^*_{\sigma }(Y)
  \]
  with grading shifted by $dim X - dim Y (mod \ 2)$.
  \end{theorem}
 
  We hope that the following remark will
   clarify some confusions  in the current literature about the role of automorphisms 
   of twists and their action on twisted  K-theroy.
   
  \begin{remark}
  \begin{enumerate}
  \item Automorphisms of the $PU(\cH)$-principal bundle $ \cP_\sigma$ 
  result in a natural action of the Piccard group  $Pic(Y)$ 
  on the twisted K-group $K^*_{\sigma }(Y)$.  Note that $Pic(Y)$ acts naturally
  on $K^*_{f^*\sigma + W_3(f)}(X)$ via the pull-back construction. Our 
  push-forward map in the  theorem is defined with a fixed choice of a
  $PU(\cH)$-principal bundle  $\cP_\sigma$ over $Y$,
   and compatible with the action of the automorphisms of $\cP_\sigma$ via
   the natural action of  $Pic (Y)$
   \[
   f_{!} ([f^*L]\cdot \xi) = [L]\cdot    f_{!} (\xi),
   \]
   for $[L]\in Pic (Y)$ and $\xi \in K^*_{f^*\sigma + W_3(f)}(X)$.  
  \item If $f$ is K-oriented, that is, $W_3(f) = 0$, our choice of  the
  Clifford bundle  $Cl( TX \oplus f^*TY)$ determines  a canonical choice of $Spin^c$  structure on
  $TX \oplus f^*TY$ with the associated first Chern class   $c_1 \in H^2(X, \ZZ)$ of the
  determinant bundle such that
  \[
  c_1 \equiv w_2(X) + f^*w_2(Y).
  \]
  Hence, we can denote $f_!$ by $f_!^{c_1}$. Then for any other choice of $Spin^c$  structure
  corresponding to $c_1 + 2c$ for  $c\in H^2(X, \ZZ)$,  the associated
  push-forward map 
  \[
  f_!^{c_1+2c}: K^*_{f^*\sigma }(X) \longrightarrow K^*_{\sigma }(Y),
  \] 
  is given by $	f_!^{c_1+2c}(\xi) = f_!^{c_1}([L_c]\cdot \xi)$, where $\xi \in K^*_{f^*\sigma}(X)$,
  and $[L_c] \in K^0(X)$ is the K-element defined by the equivalence class of the line
 bundle over $X$ with first Chern class $c\in  H^2(X, \ZZ)$.
   \end{enumerate} 
  \end{remark}

  Assume that    $dim X - dim Y  = 0 (mod \ 2)$. The proof will be divided into two parts.
        We can choose an embedding $\iota: X \to S^{2n}$, where $S^{2n}$ is an  even
   dimensional sphere. Let $g= (f, \iota): X
   \to Y \times S^{2n}$ and $\pi$ be the projection $p: Y\times
   S^{2n} \to Y$. Then $f= \pi \circ g$
   \ba\label{triangle}
   \xymatrix{ & Y\times S^{2n}\ar[dr]^{p} &\\
   X\ar[ur]^{g} \ar[rr]^{f} & & Y.
}
\na

 In the first part, we show that the
   push-forward map
  \[
   g_{!}: K_{f^*\sigma + W_3(f)}(X) \to K_{\pi^*\sigma} (Y \times
   S^{2n}),
   \]
   exists for   the smooth embedding $g$.
    Note that $f^*\sigma = g^* (\pi^*\sigma )$
   and $W_3(g)= W_3(f)$.  We will
   give a construction of the push-forward map for any smooth
   embedding in the first subsection.
In the second part, we will establish the push-forward map for any
$Spin^c$ fibration which includes the projection  $p: Y\times
   S^{2n} \to Y$. This  will be done in the second subsection.
   
As twisted K-theory is homotopy invariant, the push-forward
map is independent of the choice of the embedding $\iota: X \to S^{2n}$ for a fixed
$n$. We will show that the push-forward
map is independent of the 
dimension of the sphere $S^{2n}$ when we establish the functoriality for
push-forward maps in the third subsection.

Replacing $f: X\to Y$ by $f\times Id: X\times \RR \to Y\times \RR$, the
push-forward map of the   differentiable map $f\times Id$,
\[
(f\times Id)_!: K^*_{f^*\sigma + W_3(f)}(X\times \RR) \longrightarrow K^*_{\sigma }(Y\times \RR),
\]
induces the push-forward map for the odd degree twisted K-group
\[
f_!:  K^*_{f^*\sigma + W_3(f)}(X) \longrightarrow K^*_{\sigma }(Y).
\]

If  $dim X - dim Y  = 1 (mod \ 2)$, we can think $f: X\to Y$ as the composition of $
f\times pt: X \to Y \times S^1$ and the projection  $p: Y\times S^1\to Y$
 \[
\xymatrix{ X\ar[r]^{f\times pt } \ar[d]_{f} &Y \times S^{1} \ar[dl]^{p}\\
Y,
}
\]
then we can define the push-forward map $f_!$ as the following composition
\[\xymatrix{
f_!:  K^*_{f^*\sigma + W_3(f)}(X) \ar[r]^{\ \ (f\times pt)_!}& K^{*+1}_{p^*\sigma}(Y\times
S^1) \ar[r]^{p_!}& K^{*+1}_\sigma(Y).}
\]

\begin{corollary}
If $f: X\to Y$ is K-oriented, i.e., $W_3(f) =0$,  with $c_1\in
H^2(X, \ZZ)$ satisfying
  \[
  c_1  \equiv w_2(X) -f^*w_2(Y) \ (mod  \ 2),
  \]
  there exists a natural push-forward map
  \[
  f_{!}^{c_1}: K^*_{f^*\sigma }(X) \longrightarrow K^*_{\sigma }(Y)
  \]
  which agrees with the push-forward map for the ordinary K-theory when $\sigma$ is trivial.
  \end{corollary}

\begin{corollary} Let 
$f: X\to Y$  be a fibration with a closed oriented fiber, then the topological 
  index in twisted K-theory is given by 
  \[
  \xymatrix{
  K^*_{\pi^*f^*\sigma} (T(X/Y)) \ar[r]^{\cong } 
  &  K^*_{f^*\sigma + W_3(f)}(X) \ar[r]^{\ \ f_!} & K^*_{\sigma }(Y),
  }
  \]
   where $\pi: T(X/Y) \to X$ is the vertical tangent bundle of $f: X\to Y$. 
 \end{corollary}  

\subsection{Push-forward map for  embeddings}

Let $f: X \to Y$ be any differentiable  embedding and $\sigma
\in H^3(Y, \ZZ)$. We will show that there is a natural
homomorphism
\[
f_{!}: K_{f^*\sigma + W_3(f) }(X) \longrightarrow K_\sigma(Y).
\]
Identifying the normal bundle $V_X$ with a neighborhood of $X$ in
$Y$, so we have an embedding $\iota: V_X \to Y$. Then we have
\[
\xymatrix{ V_X\ar[r]^{\iota} \ar[d]_{\pi} & Y \\
X\ar[ur]_f.
}
\]
Denote the zero section of $V_X$  by $\underline{0}$.

\begin{lemma}  $\underline{0}^* \circ \iota^* \sigma = f^*\sigma$,
$W_3(\iota) =0$ and $W_3(\underline{0}) = W_3(f)$.
\end{lemma}
\begin{proof} $\underline{0}^* \circ \iota^* \sigma = f^*\sigma$
follows from $ \iota \circ  \underline{0} =f $. Let us calculate
$W_3(\underline{0})$:
\[\begin{array}{lll}
W_3(\underline{0}) & = & W_3(TX \oplus  \underline{0}^* T(V_X))
\\[2mm]
&=& W_3( TX \oplus \underline{0}^* (\pi^*V_X) \oplus
\underline{0}^* (\pi^*TX)) \\[2mm]
&=& W_3(TX \oplus TX \oplus V_X)
\\[2mm]
&=& W_3(V_X) = W_3(f).
\end{array}
\]
As $TV_X \oplus \iota^*TY$ admits a $Spin^c$ structure, hence,
$W_3(\iota) =0$.
\end{proof}

We define the push-forward $$\underline{0}_{!}: K_{f^*\sigma +
W_3(f) }(X) \to K_{\iota^*\sigma} (V_X) \cong K_{\pi^*\circ
f^*\sigma}(V) $$ to be  the Thom isomorphism for the oriented real
vector bundle $V_X$  equipped with  a positive definite quadratic form.

As $\iota: V_X\to Y$ is an open embedding, from the definition of
twisted K-theory, we know that there is a natural extension
map
\[
\iota_!: K_{\iota^*\sigma} (V_X) \longrightarrow K_{\sigma}(Y).
\]
This can be seen as follows:  any element in $K_{\iota^*\sigma}
(V_X) $ is represented by a
 section of $\cP_\sigma|_{V_X}\times_{PU(\cH)}
Fred^{(0)} (\cH) \to V_X$ such that away from a compact set, it is
invertible, where $\cP_\sigma$ is the principal $PU(\cH)$-bundle
over $Y$ with Dixmier-Douady class $\sigma$. By homotopy, such a
section can be chosen so that away from a compact set, it is the
identity operator, hence, naturally defines a section of
$\cP_\sigma\times_{PU(\cH)} Fred^{(0)} (\cH)\to Y$.

We define the push-forward map for $f: X \to Y$ to be
\ba\label{push:close}
 f_! = \iota_! \circ \underline{0}_{!}.
 \na
 Applying the homotopy invariance of twisted K-theory, we know
 that $f_!$ doesn't depend on the choice of the identification of
 the normal bundle $V_X$ with a tubular neighborhood $X$ in $Y$.

\subsection{Push-forward map for $Spin^c$ fibrations}

Given  a smooth fibration  $\pi: X \to Y$ over a closed smooth manifold $Y$,
whose fibers are diffeomorphic to an even dimensional
compact, oriented, closed $Spin^c$ manifold $M$,
we will define a push-forward map
on twisted $K$-theory using a twisted family of
Fredholm operators coupled to Dirac operators along the fibre of $\pi$. 

In fact, 
from  the algebraic definition of the twisted K-theory, one sees that the
push-forward map 
\[
\pi_!: K^0_{\pi^*\sigma} (X) \longrightarrow  K^0_{\sigma}(Y),
\]
is defined via the Kasparov product with the longitudinal
Dirac element, which is the element
in $KK(C_0(X), C_0(Y))$  defined by the fiberwise Dirac
operator (Cf. \cite{CS}).  Let us illustrate this idea using the following projection 
\[
p: Y\times  S^{2n} \longrightarrow Y, 
\]
which is enough  for our construction of the push-forward map (Cf. Diagram \ref{triangle}).

 Recall the algebraic definition of the twisted K-theory (Cf. \cite{AS} \cite{Ros}). 
  Let $\cH$ be a infinite dimensional, complex and separable Hilbert space, and 
   $\cK(\cH)$ be the $C^*$-algebra of  compact   operators on $\cH$.
 Given a class $\sigma\in H^3(Y, \ZZ)$ on a locally compact Hausdorff
  space $Y$, we can associate a principal $PU(\cH)$-bundle
  $\cP_\sigma$ over $Y$ with its Dixmier-Douady class given by $\sigma$. Applying the
  natural identification $PU(\cH) \cong Aut(\cK(\cH))$ (the *-automorphism group
  of $\cK(\cH)$), one can  has a bundle of compact operators 
  \[
  \cP_\sigma\times_{PU(\cH)} \cK(\cH),
  \]
  of which the space of  continuous sections vanishing at infinity
   is denoted by $\AA_\sigma (Y)$. 
 
 \begin{definition} The algebraic definition of  twisted K-theory 
  is  given by  the K-theory of the corresponding
  continuous trace $C^*$-algebra $\AA_\sigma (Y)$
  \[
  K^i_\sigma(Y) := K^i_{\cP_\sigma} (Y) = KK^i(\CC, \AA_\sigma(Y) ),
  \]
  for $i=0, 1$. 
  \end{definition}

  \begin{proposition} For the natural projection $p: Y\times  S^{2n} \longrightarrow Y$,
  the push-forward map 
  \[p_!:  K^0_{p^*\sigma }( Y\times S^{2n}) \longrightarrow K^0_\sigma(Y)
  \]
  is given by the Kasparov product with the Dirac  element 
  $[\delta_{S^{2n}}] \in KK(C(S^{2n}), \CC)$. 
  \end{proposition}
  \begin{proof} 
Denote by $p^*\sigma \in H^3(Y\times  S^{2n}, \ZZ)$ the pull-back of $\sigma$ via 
$p:Y\times  S^{2n} \to Y$, 
then we have
\[
\AA_{p^*\sigma} (Y \times S^{2n})\cong C(S^{2n} ) \otimes \AA_\sigma(Y).
\]
It is well-known that the Dirac operator on $S^{2n}$ with  standard metric defines
a canonical  analytic K-cycle:
\[
[\delta_{S^{2n}}] \in KK(C(S^{2n}), \CC).
\]
Applying the Kasparov product,
\[
KK\bigl(\CC, C(S^{2n}) \otimes \AA_\sigma(Y)\bigr) \otimes KK(C(S^{2n}), \CC)
\longrightarrow  KK(\CC, \AA_\sigma(Y)),
\]
we get the push-forward map associated to $p: Y\times   S^{2n} \longrightarrow Y$,
\[
p_! = \otimes  [\delta_{S^{2n}}]: K^0_{p^*\sigma }( Y\times S^{2n}) \cong 
KK(\CC, \AA_{p^*\sigma} (Y \times S^{2n}) ) \longrightarrow  K^0_\sigma(Y).
\]
\end{proof}

As we use the model of twisted Fredholm operators
to define the twisted K-theory, we now provide  a geometric proof of the existence
of the push-forward map, which is merely the obvious translation of
the Kasparov product to the differential geometry set-up. 
For simplicity, 
we also assume that $X$ and $Y$ are compact manifolds.

Let $g^{X/Y}$ be a metric on the relative tangent bundle $T(X/Y)$
(the vertical tangent subbundle of $TX$). The spinor bundle associated to
$(T(X/Y), g^{X/Y})$ is denoted by
$$S_{X/Y} = S_{X/Y}^+ \oplus S_{X/Y}^-.$$
 Let $T^HX$ be a smooth vector subbundle of $TX$, (a horizontal vector
 subbundle of $TX$),  such that
 \ba\label{horizontal}
 TX = T^HX \oplus T(X/Y).
 \na
Denote by $P^v$ the projection of $TX$ onto the vertical tangent bundle
 under the decomposition (\ref{horizontal}).
As shown in Theorem 1.9 of \cite{Bis1}, $(T^HX, g^{X/Y})$ determines a
canonical Euclidean
 connection $\nabla^{X/Y}$ as follows. Choose a metric $g^X$ on $TX$ such that
 $T^HX$ is orthogonal to $T(X/Y)$ and such that $g^{X/Y}$ is
the restriction of
 $g^X$ to $T(X/Y)$. Then
 \[
 \nabla^{X/Y} = P^v \circ \nabla^X
 \]
 where $\nabla^X$ is the Levi-Civita connection on $(TX, g^X)$.

The Clifford multiplication of $T^*(X/Y)$ on $S_{X/Y}$ and the Euclidean
connection $\nabla^{X/Y}$ defines a fiberwise Dirac operator acting on
$C^\infty (X_y,  S_{X/Y})$, which can be completed to a Hilbert space,
denoted by $\cH_{y}$.
For any $y\in Y$, let $\Dirac_y$ be the Dirac operator along the fiber $\pi^{-1}(y)$.
Then $\{\Dirac_y\}_{y\in Y}$
is a  family of
elliptic  operators over $Y$ acting on an infinite dimensional
Hilbert bundle $\cH_{X/Y}$, whose fiber at $y$ is $\cH_y$.
Note that $\cH_{X/Y}$ is a $\ZZ_2$-graded Hilbert bundle
over $Y$ for the $Spin^c$ fibration with even dimensional fibres.

 The partial isometry part in the polar decomposition of $\Dirac_y$ is not a continuous family of
  bounded Fredholm operators parametrized by $y\in Y$. Modulo compact operators,  there
  is  a continuous family of bounded Fredholm operators parametrized by $y\in Y$,
 denoted by 
 \ba\label{bded} V_y = \Dirac_y ( \Dirac_y^* \Dirac_y +1 )^{-\frac 12},
 \na
 under the usual assumption that $\Dirac_y$ is a smooth family of Dirac operators parametrized
 by $Y$.

Take a  trivialization
of $P_\sigma$ for a cover $\{U_\alpha\}$ of $Y$:
\[
P_\sigma |_{U_\alpha} \cong U_\alpha \times \cP(\hat{\cH}),
\]
with the transition functions given by $g_{\alpha\beta}: U_{\alpha\beta} \to PU(\hat{\cH}).$
We assume that the fibration $\pi: X \to Y$ is also trivialized over $U_\alpha$, i.e.,
\[
\pi^{-1} (U_\alpha ) \cong U_\alpha \times M.
\]
Then there is a trivialization of $P_{\pi^*\sigma}= \pi^*  P_\sigma$ with respect to
the cover  $\{ \pi^{-1} (U_\alpha )\}$ of $X$ and 
the transition functions given by $\pi^*g_{\alpha\beta} =g_{\alpha\beta}\circ \pi$. 

Given   an element $T \in K^0_{P_{\pi^*\sigma}}(X)$, we can  represent it by
a  section of $$P_{\pi^*\sigma}\times_{PU(\hat{\cH})} Fred^{(0)}(\hat{\cH}),$$
written as 
\[
T_\alpha: \pi^{-1} (U_\alpha ) \longrightarrow Fred^{(0)}(\hat{\cH})
\]
satisfying $ T_\beta = (\pi^*g_{\alpha\beta}) T_\alpha  (\pi^*g_{\alpha\beta}) ^{-1}.$
  We can couple the fiberwise Fredholm  operators 
  $\{V_y\}$ in (\ref{bded})  with $T_\alpha$ to get a section of
  $$
  Fred^{(0)}(P_{\sigma}) = 
  P_{\sigma}\times_{PU(\hat{\cH})}Fred^{(0)}(\hat{\cH}\otimes \cH_{X/Y} ), 
 $$
  by graded tensor product
  \[
  \pi_!(T_\alpha)= V\otimes I + I \otimes  T_\alpha: U_\alpha \longrightarrow
Fred^{(0)}(\hat{\cH}\otimes \cH_{X/Y}),
  \]
 as these local   sections    $\{ \pi_!(T_\alpha) \}$    satisfy  
  $$ 
  \pi_!(T_\beta) =(g_{\alpha\beta}\times Id) \pi_!(T_\alpha)(g_{\alpha\beta}\otimes Id)^{-1}.
  $$ 
 
 Fix a graded isomorphism $\phi: \hat{\cH} \otimes  \cH_{X/Y} \to  \hat{\cH}$. We have
 a principal $PU(\hat{\cH})$-bundle $\tilde{P}_{\sigma}$  whose
 transition functions with respect to the cover $\{U_\alpha\}$ of $Y$
 are given by 
 \[
\tilde{g}_{\alpha\beta} = \phi \circ (g_{\alpha\beta}\times Id) \circ \phi^{-1}: U_{\alpha\beta} 
 \longrightarrow PU(\hat{\cH}).
 \]
 Hence, $\tilde {P}_{\sigma} \cong {P}_{\sigma}$.  Then the above
 local section $\{ \pi_!(T_\alpha) \}$ induces a section of 
 $$Fred^{(0)}( \tilde{P}_{\sigma}) =\tilde{P}_{\sigma} \times_{ PU(\hat{\cH})} 
  Fred^{(0)}(\hat{\cH}).$$
 Notice that a homotopy equivalent section of   $Fred^{(0)}( P_{\pi^*\sigma})$,
  under the above procedure, produces  a homotopy equivalent section
  of  $Fred^{(0)}( \tilde{P}_{\sigma})$. 
  
  A different choice of the graded isomorphism 
  $\phi: \hat{\cH} \otimes  \cH_{X/Y} \to  \hat{\cH}$ gives
  rise to an isomorphism on the resulting twisted K-group. Note that
  $Aut (P_{\sigma}) \cong H^2(Y, \ZZ)$ acts on $K^0_{P_\sigma}(Y)$ by natural
   isomorphisms.   This finishes the definition of  the push-forward map, 
   up to the natural actions of
 $H^2(X, \ZZ)$  and  $H^2(Y, \ZZ)$, 
  $$\pi_!: K^0_{ \pi^*\sigma } (X) \to  K^0_{ \sigma}(Y).$$

  The  push-forward map on the odd degree $K$-group
  $$ \pi_!: K^1_{ \pi^*\sigma } (X) \to  K^1_{ \sigma}(Y) $$
  can also be defined using a based $S^1$-family of local twisted sections.

\begin{remark} For a $Spin^c$ fibration $\pi: X \to Y$ with odd dimensional fibers,
with the help of isomorphisms (\ref{0:1}) and (\ref{1:0}), we obtain the corresponding
push-forward maps:
\ba\label{push:odd}
\begin{array}{c}
\pi_!:  K^1_{P_{\pi^*\sigma}} (X) \to  K^0_{P_\sigma}(Y),\\[2mm]
\pi_!: K^0_{P_{\pi^*\sigma}} (X) \to  K^1_{P_\sigma}(Y).
\end{array}
\na
 \end{remark}

\subsection{Functoriality for push-forward maps}

In this subsection, we will establish the wrong way functoriality 
\[
(g\circ f)_! = g_! \circ f_!:  K^*_{(g\circ f)^*\sigma + W_3(g\circ f)}(X) 
\to  K^*_{g^*\sigma + W_3(g)  } (Y) \to K^*_\sigma (Z) 
\]
for two differentiable maps $f: X\to Y$, $ g: Y\to Z$ and $\sigma\in H^3(Z, \ZZ)$. 

We first prove  the wrong way functoriality for two embeddings.

\begin{lemma}\label{func:emb} Let  $f: X\to Y$ and  $ g: Y\to Z$ be two smooth embeddings, then 
$(g\circ f)_! = g_! \circ f_!$.
\end{lemma}
\begin{proof} For a smooth embedding, the corresponding push-forward map is given
by the composition of Thom isomorphism and the natural extension for
open embeddings. Denote by $V_{X\subset Y} $ and $V_{Y\subset Z}$  the normal
bundles over  $X$   and $Y$,   together with fixed
identifications to  the tubular neighbourhoods of $f(X)$ and $g(Y)$. Denote by $\underline{0}_X$
and  $\underline{0}_Y$ the corresponding zero sections of $V_{X\subset Y} $ and
 $V_{Y\subset Z}$ respectively, then we have the following commutative diagram
\[
 \xymatrix{  && i_1^* ( V_{Y\subset Z})  \ar[dr]^{ i_3}&& \\
   & V_{X\subset Y}\ar[dr]^{i_1}\ar[ur]^{\underline{0} } &  & 
 V_{Y\subset Z}  \ar[dr]^{i_2}& \\
   X\ar[ur]^{\underline{0}_X } \ar[rr]^{f} & & Y\ar[ur]^{ \underline{0}_Y}  
    \ar[rr]^{g}  && Z,}
\]
    where 
    $i_1^*  (V_{Y\subset Z} )$ is  the restriction of $V_{Y\subset Z} $ 
    to the open tubular neighbourhood $V_{X\subset Y}$ , and $i_1$, $i_2$ and $i_3$
    are open embeddings.
    As $i_1^* ( V_{Y\subset Z} )$ is also a normal bundle over   $X$, with a fixed
    identification to the tubular neighbourhood of $(g\circ f) (X)$ in $Z$, then 
$(g\circ f)_! = g_! \circ f_!$ follows from the definition of the push-forward map.
\end{proof}

Now we show functoriality for general smooth maps $f: X\to Y$ and
 $ g: Y\to Z$. Choosing two embeddings $\iota_1: X \to S^{2m}$ and
$\iota_2: Y \to S^{2n}$, we have the following commutative diagram:
   \ba\label{triangle2}
   \xymatrix{  && Z\times S^{2n} \times S^{2m} \ar[dr]^{\pi_4} 
   \ar@/_2pc/@{-->}[ddd]^{\pi_5}&& \\
   & Y\times S^{2m}\ar[dr]^{\pi_1}\ar[ur]^{ g\times \iota_2\times Id} 
   &  & Z\times  S^{2n}\ar[dr]^{\pi_2}& \\
   X\ar[ur]^{f \times \iota_1} \ar[drr]_{(g\circ f) \times \iota_1}\ar[rr]^{f} &
    & Y\ar[ur]^{g \times \iota_2}  
    \ar[rr]^{g}  && Z \\
  && Z\times S^{2m}\ar[urr]_{\pi_3} && ,
}
\na
where $\pi_5$ is the natural projection  
$  Z\times S^{2n} \times S^{2m} \to Z\times S^{2m}$. 
 Then from the definition of the push-forward map, we know
 that  
 \[\begin{array}{c}
 g_! \circ f_! = (\pi_2)_! \circ (g\times \iota_2)_! \circ (\pi_1)_! \circ (f \times \iota_1)_!\\
 (g\circ f)_! =   (\pi_3)_! \circ\bigl( (g\circ f) \times \iota_1\bigr)_!.
 \end{array}
  \]
  The proof of $(g\circ f)_! = g_! \circ f_!$ 
  is given in the following  sequence of identities.
  
  \begin{lemma} \label{square} $(g\times \iota_2)_! \circ (\pi_1)_! =  
   (\pi_4)_! \circ (g\times \iota_2 \times Id)$.
  \end{lemma}
  \begin{proof}
  The proof of this lemma  follows directly from the definition of the
  push-forward map for embeddings and the following
  commutative diagram:
  \ba\label{square:triangle}
  \xymatrix{  & V_Y \times S^{2m}\ar@{-->}[dd] \ar[dr]  &\\
  Y\times S^{2m}\ar[dd]^{\pi_1}  \ar[ur]^{\underline{0}} \ar[rr]_{\qquad g\times \iota_2\times Id} &  &
  Z\times S^{2n}\times S^{2m}  \ar[dd]^{\pi_4} \\
  & V_Y \ar@{-->}[dr] &\\
  Y  \ar[rr]^{g\times \iota_2} \ar@{-->}[ur]^{\underline{0}}  & & Z \times S^{2n}.}
  \na 
  where the normal bundle $V_Y$ is identified with the tubular neighbourhood of $Y$ in
  $ Z\times S^{2n}$, with zero section $\underline{0}:  Y\to V_Y$, 
  and the normal bundle $V_Y\times S^{2m}$ is identified with the tubular 
  neighbourhood of $Y\times S^{2m}$ in
  $ Z \times S^{2m}\times S^{2n}$ .
  \end{proof}
  
  The  following  commutative diagram with the obvious projections
  \[
  \xymatrix{   
   Z\times S^{2n}\times S^{2m}  \ar[r]^{\pi_4} \ar[d]^{\pi_5} 
    &  Z\times S^{2n}    \ar[d]^{\pi_2} \\
   Z\times S^{2m}  \ar[r]^{\pi_3} & Z
   }
   \]
   implies that $(\pi_2)_! \circ (\pi_4)_! = (\pi_3)_! \circ (\pi_5)_!$, as
   the Dirac element on $S^{2m}\times S^{2n}$ is given by
   the external product of the  Dirac elements on $S^{2m}$ and on $S^{2n}$.
   
   Putting together the above identities, we know that
     $(g\circ f)_! = g_! \circ f_!$ follows,  if we can show
     that
     \ba\label{id}
  (\pi_3)_!\circ  \bigl( (g\circ f)\times \iota_1\bigr)_! = 
 (\pi_3)_!\circ  (\pi_5)_! \circ (g\times \iota_2\times Id)_!\circ  (f\times \iota_1)_!. 
     \na
  Choose an isometric embedding $\iota: S^{2n}\times S^{2m} \to S^{2N}$ where these
  spheres are equipped with the natural round metrics.
    Consider the following  commutative diagram, where maps in the top row
    are all smooth embeddings,  
     \[\xymatrix{  X \ar[r]^{f\times \iota_1\ }  
     \ar[drr]_{(g\circ f)\times \iota_1} &  Y \times S^{2m}
      \ar[r]^{g\times \iota_2\times Id \quad }
     &  Z\times S^{2n}\times S^{2m} \ar[d]^{\pi_5}\ar[r]^{Id \times \iota}   &
     Z \times S^{2N} \ar[d]^{\pi_6} \\
  &&   Z\times S^{2m}  \ar[r]^{\pi_3} & Z.
     }
     \]
     Notice that the right square in the above diagram tells us 
     that 
     \[
     (\pi_6)_! \circ (Id \times \iota )_! =  (\pi_3)_! \circ (\pi_5)_!.
     \]
     Here we use  functoriality in  topological index theory
     \[
     \xymatrix{ 
     S^{2m}\times S^{2n} \ar[r]^{\iota}\ar[d] & S^{2N} \ar[d]\\
     S^{2n}\ar[r] & pt,
     }\]
   which reads,   in terms of Dirac elements on the round spheres, 
   \[
   [\delta_{S^{2m}\times S^{2n}} ] = \iota !\otimes [\delta_{S^{2N}}],
   \]
   under the Kasparov product (Cf. \cite{CS}).  Here $\iota !$ is the canonical
   element in KK-thory $$KK(C(S^{2m}\times S^{2n}), C(S^{2N}))$$
    defined in \cite{CS} by Connes and Skandalis. 
   Hence, the required equation  (\ref{id}) becomes
     \[
      (\pi_3)_!\circ  \bigl((g\circ f)\times \iota_1\bigr)_! =
    (\pi_6)_! \circ  (Id \times \iota)_! 
    \circ (g\times \iota_2\times Id)_!\circ  (f\times \iota_1)_!.
    \]
Apply  functoriality for embedding maps (Cf. Lemma \ref{func:emb}), so
     we know that 
     \[
     (Id \times \iota)_! \circ (g\times \iota_2\times Id)_!\circ  (f\times \iota_1)_!
     = \bigl((g\circ f) \times \tilde{\iota}\bigr)_!,
     \]
     where $(g\circ f) \times \tilde{\iota}: X \to Z\times S^{2N}$ is the induced
     embedding from $\iota_1$, $\iota_2$ and $\iota$.

     It remains to show that 
     \[
     (\pi_3)_!\circ  \bigl((g\circ f)\times \iota_1\bigr)_! =
     (\pi_6)_! \circ \bigl((g\circ f) \times \tilde{\iota}\bigr)_!. 
     \]
     Equivalently, the push-forward map $(g\circ f)_!$
     is independent of the choice of sphere in its definition. We establish this fact for
     a general smooth map $f: X \to Z$, with a slight change of notations
     for simplicity.
     For two different embeddings $\iota_1: X\to  S^{2m}$ and
          $\iota_2: X\to  S^{2n}$, we can choose smooth embeddings 
          $ \iota': S^{2m} \to S^{2N}$ and $ \iota'': S^{2n} \to S^{2N}$. 
          Then the proof follows from the following
          commutative diagram:
          \[         
    \xymatrix{ & Y\times S^{2m}\ar[d] \ar[dr]^{Id\times \iota'}  &\\
    X\ar[r]^{f}\ar[ur]^{f\times \iota_1} \ar[dr]_{f\times\iota_2} & Y  & Y\times S^{2N}.\ar[l]\\
    &  Y\times S^{2n}\ar[u] \ar[ur]_{Id\times \iota''} &}
    \]
  Here we consecutively  apply  functoriality in  topological index 
   theory for the following diagrams: 
   \[
\xymatrix{  S^{2m}\ar[d] \ar[dr]^{\iota'} &\\ pt & S^{2N}\ar[l] 
  } \qquad  \qquad \qquad 
  \xymatrix{ pt &  S^{2N} \ar[l]\\
  S^{2n}\ar[u]\ar[ur]_{\iota''}&, }
  \]
i.e., $ [\delta_{S^{2m}}] = (\iota') ! \otimes [\delta_{S^{2N}}]$, 
$ [\delta_{S^{2n}}] = (\iota'') ! \otimes [\delta_{S^{2N}}]$, the  
functoriality  for embedding maps (Lemma \ref{func:emb}) 
\[
(Id \times \iota') \circ (f\times \iota_1) = f\times (\iota'\circ \iota_1),
\qquad 
(Id \times \iota'') \circ (f\times \iota_1) = f\times (\iota''\circ \iota_2)
\]
and the homotopy invariance of   twisted K-theory for two embeddings 
$\iota'\circ \iota_1$ and $\iota''\circ \iota_2$.

\section{Some applications}

\subsection{Index of Dirac operators for Clifford modules}

According to the Atiyah-Singer index theory for a Dirac operator
on a Clifford module $E$ over a compact oriented even dimensional
manifold $X$, see Theorem 4.3 \cite{BGV}, the index is given by
\ba\label{AS:index}
 \disp{\int_X}\hat{A}(X)Ch(E/S),
\na
 where $Ch(E/S)$ is the relative Chern character form defined
by a Clifford connection on $E$.

Let $Cl(X)$ be the bundle of Clifford algebra associated to the
tangent bundle equipped with the Riemannian metric. From Lemma
\ref{Cliff=gb}, we know that the category of $Cl(X)$-modules is
isomorphic to the category $\cE_{bg}(X, \Gamma_{W_3(X)})$ of
$\Gamma_{W_3(X)}$-modules.
We can express the Atiyah-Singer index
(\ref{AS:index}) as the topological index:
\[
t-index: K_{W_3(X)}(X) \longrightarrow \ZZ.
\]
By the Thom isomorphism for the tangent bundle equipped with a
Riemannian metric (Cf. Theorem \ref{Thom1} ), we have
\[ K_{W_3(X)}(X) \cong K(TX).\]
 The topological index is given by
\ba\label{topo:index}
 K(TX) \longrightarrow \ZZ.
\na
 Recall the definition of (\ref{topo:index}) for a compact
 manifold $X$ by choose an embedding $\iota: X \to \RR^m$, then
 the tangent map, also denoted by
 $$
 \iota: TX \longrightarrow  T\RR^m \cong 
 \CC^m,
 $$
 is K-oriented . Then the corresponding push forward map
 in ordinary K-theory and Bott periodicity give rise to
 \[
\iota_!: K(TX) \longrightarrow K(T\RR^m) \cong K(\CC^m) \cong \ZZ.
\]
This proves that the Atiyah-Singer index (\ref{AS:index}) is the
push-forward map in twisted K-theory for $f: X \to pt$.

In \cite{MurSin}, generalized Dirac operators are introduced for
any even dimensional  Riemannian manifold $M$ in terms of bundle
gerbe modules for the lifting bundle gerbe $\Gamma_{W_3(X)}$.
 The index of Murray-Singer's generalized
Dirac operators defines a map:
\[
Ind: \cE_{bg}(X, \Gamma_{W_3(X)}) \longrightarrow \ZZ,
\]
which agrees with  Atiyah-Singer index (\ref{AS:index}) under the
category equivalence
$$\cE_{bg}(X, \Gamma_{W_3(X)} )\cong \cE^{TX}(X).$$
Descending to the twisted K-group $K_{ W_3(X)}(X)$, the
 Murray-Singer index is the push-forward map in twisted K-theory for  $f: X\to
 pt$.

\subsection{D-brane charges}

In type IIA/IIB superstring theory with topologically trivial
$B$-field, a $D$-brane (see \cite{Wit1}) in an oriented
10-dimensional $Spin$ manifold $X$, is a $Spin^c$ submanifold
$\iota: Q \to X$, together with a Chan-Paton bundle and a
superconnection, defined by an element $\xi\in K(Q)$. Using the
push-forward map in ordinary K-theory, $$\iota_!: K(Q)
\longrightarrow K^*(X)$$ we see that such $D$-branes are
classified by elements in $K(X)$, and the Ramond-Ramond charge of
$(Q, \xi)$ is given by
\[
Ch(\iota_!\xi) \sqrt{\hat{A}(X)} \in H^*(X, \RR).
\]

When the $B$-field is non-trivial and its  characteristic class
$\sigma$ lies in $ H^3(M, \ZZ)$ with the curvature denoted by  $H$,
there is a constraint.
A submanifold $\iota: Q \to X$, has associated to it a well defined action when it satisfies
the Freed-Witten anomaly cancellation condition
(Cf. \cite{FreWit}):
\[
\iota^*\sigma + W_3(Q)=0,
\]
in $H^3(Q, \ZZ)$. It was proposed in
\cite{Wit2}\cite{Kap}\cite{BouMat} that $D$-brane charges should
be classified by the twisted K-group $K^*_\sigma(X)$, yet a
rigorous formulation of such $D$-brane charges hasn't been found.

Now as an application of our push-forward map in twisted K-theory,
we can associate a canonical element in $K_\sigma(X)$, which can
be called the $D$-brane charge of the underlying $D$-branes.  As
$X$ is spin, $W_3(Q)$ agrees with $W_3(V_Q)$, the third
Stiefel-Whitney class of a normal bundle of $Q$ in $X$. Hence,
$W_3(Q) = W_3(\iota)$.  For a submanifold $\iota:  Q \to X$
satisfying the Freed-Witten anomaly cancellation condition $
\iota^*\sigma + W_3(Q) =0, $ in $H^3(Q, \ZZ)$.  We apply the
push-forward map for $\iota$ to get
\[
\iota_!: K(Q) \cong K_{\iota^*\sigma + W_3(Q)} (Q) \longrightarrow
K^*_\sigma (X),
\]
hence, for any $D$-brane wrapping on $Q$ determined by an element
$\xi\in K(Q)$, we can define its charge as
\[
\iota_!(\xi) \in K^*_\sigma(X).
\]


\begin{thebibliography}{9999}


\bibitem{AH}  M. Atiyah, F. Hirzebruch {\sl Vector bundles and homogeneous spaces},
Proceedings of Symposium in Pure Mathematics, Vol. 3, 7-38.  Am. Math. Soc. 1961.

\bibitem{AS} M. Atiyah, G. Segal {\sl Twisted $K$-theory}, preprint.



\bibitem{BEM} P. Bouwknegt,  J. Evslin, V.  Mathai  
   {\sl T-Duality: Topology Change from H-flux.  }   Commun.Math.Phys. 249 (2004) 383-415



\bibitem{BGV} N. Berline, E. Getzler, M. Vergne  {\sl Heat kernels and Dirac operators},
Springer-Verlag. Berlin, 1992.

\bibitem{Bis1} J.M. Bismut  {\em  Local index theory, eta invariants and
holomorphic torsion: a survey.}  Surveys in differential geometry, Vol. III,
1--76, Int. Press, Boston, MA, 1998.


\bibitem{BCMMS} P. Bouwknegt, A. Carey, V. Mathai, M.  Murray, D.
 Stevenson  {\em Twisted $K$-theory and
 $K$-theory of bundle gerbes.}  Comm. Math. Phys., Vol. 228    no. 1, 17--45, 2002.

\bibitem{BouMat} P. Bouwknegt,  V. Mathai  {\it $D$-branes, $B$-fields and
twisted $K$-theory}, J. High Energy Phys. {\bf 03} (2000) 007,
hep-th/0002023.

\bibitem{CJMSW} A. Carey, S. Johnson, M. Murray, D. Stevenson, B.L. Wang,
{\sl   Bundle gerbes for Chern-Simons and Wess-Zumino-Witten theories}, 
Comm.   Math. Phys., Vol. 259  No. 3 577--613, 2005.

\bibitem{CW1} A. Carey, B.L. Wang  {\sl On the relationship of
 gerbes to the odd families index theorem,}
    preprint, math.DG/0407243. 
    
    
\bibitem{CW2} A. Carey, B.L. Wang  {\sl   Fusion of Symmetric D-branes and Verlinde ring,} 
preprint,  math.ph/0505040. 
 
 
 \bibitem{CS} A. Connes, G. Skandalis  {\sl The longitudinal index theorem for foliations.}
   Publ. Res. Inst. Math. Sci.  20  (1984),  no. 6, 1139--1183. 

 \bibitem{DK} P. Donavan, M. Karoubi   {\sl Graded Brauer groups and K-theory with lcoal
 coefficients}, Publ. Math. de IHES, Vol 38, 5-25, 1970.

\bibitem{FHT} D. Freed, M.J. Hopkins and C. Teleman {\sl
Twisted K-theory and Loop Group Representations, }  math.AT/0312155,  preprint.


\bibitem{FreWit}
D.~Freed,  E.~Witten {\it Anomalies in String Theory with
D-Branes}, hep-th/9907189.

\bibitem{Kap}
A. Kapustin  {\it $D$-branes in a topologically non-trivial
$B$-field}, Adv.\ Theor.\ Math.\ Phys.\ {\bf 4} (2001) 127,
hep-th/9909089.

 \bibitem{Kar0} M. Karoubi   {\sl   K-theory, an introduction.}
  Grundlehren der math. Wiss. N° 226. Springer Verlag (1978).

 \bibitem{LM} H. B. Lawson, M-L Michelsohn {\sl Spin Geometry}. Princeton University
 Press, 1989.

\bibitem{MMSin} V. Mathai, R. Melrose, I. Singer {\sl The index 
of projective families of elliptic operators,}  Geom. Topol. 9(2005) 341-373.


 \bibitem{MatSte} V. Mathai, D. Stevenson
 {\sl Chern character in twisted K-theory: equivariant and holomorphic cases},
 Commun.Math.Phys. 228 (2002) 17-49.



\bibitem{Mur}
M. K. Murray  {\sl Bundle gerbes,}
J. London Math. Soc. (2) {\bf 54}
(1996), no.~2, 403--416.


\bibitem{MurSin} M. K. Murray, M. A. Singer {\sl Gerbes, Clifford modules and the
 index theorem. }
  Ann. Global Anal. Geom.  26  (2004),  no. 4, 355--367.


  \bibitem{Par1} E. Parker {\sl The Brauer group of 
  Graded continuous trace $C\sp *$-algebras}. Trans. of Amer. Math. Soc., Vol 308, No. 1, 1988,
  115-132.
  
  \bibitem{Par2} E. Parker {\sl 
  Graded continuous trace $C\sp *$-algebras and duality. }
Operator algebras and topology (Craiova, 1989), 130--145,
Pitman Res. Notes Math. Ser., 270.


 \bibitem{PS} A. Pressley, G. Segal {\sl Loop groups},
    Oxford University Press, Oxford,  1988.

 \bibitem{Ros} J. Rosenberg {\sl Continuous-trace algebras from
 the bundle-theoretic point of view.}
    Jour. of Australian Math. Soc. Vol A 47, 368-381, 1989.


\bibitem{Seg} G. Segal  {\sl Equivariant K-theory}
  Publications Math. de IHES, Vol. 34, 129-151, 1968..

\bibitem{TXC} J. Tu, P. Xu, C. Laurent-Gengoux {\sl Twisted K-theory of differentiable 
stacks}, Ann. Sci. Ecole Norm. Sup. (4) Vol 37, 841-910, 2004.

\bibitem{Wit1} E. Witten  {\sl  D-branes and K-theory}
 J. High Energy Phys.  1998,  no. 12, Paper 19,

\bibitem{Wit2}
E. Witten  {\it Overview of $K$-theory applied to strings}, Int.
J. Mod. Phys. {\bf A16} (2001) 693, hep-th/0007175.


  \end{thebibliography}
  \end{document}